\documentclass[11pt]{article}
\usepackage{graphicx}
\usepackage{amsmath,amsfonts,amsthm}
\usepackage{amssymb,latexsym}
\usepackage{fixmath}
\usepackage{mathrsfs,amsbsy}
\usepackage{dsfont}
\usepackage[normalem]{ulem}
\usepackage{caption}
\usepackage{subfigure}
\usepackage{epstopdf}
\usepackage{algorithm}
\usepackage{algorithmic}
\usepackage[algo2e,ruled,vlined,linesnumbered, nofillcomment]{algorithm2e}

\usepackage{dcounter} % added <<<<<<<<<<<<<<<<<
\countstyle{section}%
\DeclareDynamicCounter{algocf}%

\usepackage{eucal}
\usepackage{colortbl}
\usepackage{multirow}
\usepackage[bottom]{footmisc}

\def\bbR{\mathbb{R}}

\def\wtd{\widetilde}

\usepackage{accents}
\newcommand\munderbar[1]{%
	\underaccent{\bar}{#1}}

\def\cR{{\cal R}}

\DeclareMathOperator{\T}{T}

\DeclareMathOperator*{\argmin}{argmin}

\DeclareMathOperator{\subspan}{span}

\newtheorem{remark}{Remark}[section]

\newtheorem{theorem}{Theorem}[section]

\newtheorem{corollary}{Corollary}[section]

\numberwithin{equation}{section}
\numberwithin{figure}{section}
\numberwithin{table}{section}

\title{Flexible Modified LSMR for Least Squares Problems}
\author{Mei Yang\thanks{Division of Data Science, College of Science, University of  Texas at Arlington, Arlington, P.O. Box 19047, TX 76019, USA.
		Email: {\tt mei.yang@uta.edu}.
	}
	\and
	Gul Karaduman \thanks{
		Department of Mathematics, Karamanoglu Mehmetbey University, 70100, Karaman, Turkey.  Email:{\tt gulk@bu.edu}.}
	\and
	Ren-Cang Li  \thanks{
		Department of Mathematics,
		University of Texas at Arlington,
		P.O. Box 19408,
		Arlington, TX 76019-0408, USA.
		Supported in part by NSF DMS-2009689.
		Email: {\tt rcli@uta.edu}.
	}
}

\date{\today
}

\begin{document}
	
	\maketitle
	
	\begin{abstract}
		LSMR  is a widely recognized method for solving least squares problems via the double QR decomposition.
		Various preconditioning techniques have been explored to improve its efficiency.
		One issue that arises when implementing these preconditioning techniques is the need to solve two linear systems per iterative step.
		In this paper, to tackle this issue, among others,
		a modified LSMR method (MLSMR), in which only one linear system per iterative step needs to be solved instead of two, is
		introduced, and then it is integrated with the idea of flexible GMRES to yield a flexible MLSMR method (FMLSMR).
		%It constructs flexible preconditioners based on solving one linear equation.
		%Its application to regularized least squares problems, commonly used in regularization models, are discussed.
		Numerical examples are presented to demonstrate the efficiency of  the proposed FMLSMR method.
	\end{abstract}
	
	\smallskip
	{\bf Key words.} LSMR, flexible GMRES, FMLSMR, preconditioner, least squares problem, linear system
	\smallskip
	
	{\bf AMS classification.} 65F10
	
	\section{Introduction}\label{sec:intro}
	%	\section{Preliminaries and Notations}
	In this paper, we present a numerical method called the flexible modified LSMR method (FMLSMR),  which is built upon our previous work \cite {Yang:2018},  for solving the least square problem
	
	\begin{equation}\label{eq:LSP}
		\min_x\left\| Ax-b\right\|_2,
	\end{equation}
	where $A\in \mathbb{R}^{m\times n}$, $b\in\mathbb{R}^m$, and either $m\ge n$ or $m< n$ is allowed.
	FMLSMR  improves of the well-known LSMR method \cite{fosa:2011} which is based on the Golub-Kahan bidiagonalization.
	LSMR seeks the best approximate solution in the Krylov subspace  $\mathcal{K}_k\left(A^{\T}A,A^{\T}b\right)$ that minimizes
	$\|A^{\T}r\|_2$, while LSQR \cite{pasa:1982} seeks to minimize $\|r\|_2$ in the same Krylov subspace, where $r=b-Ax$. LSMR has been used for various problems, including saddle point problems, as demonstrated in recent publications  \cite{Gul:2022,Gul:2023}.
	
	To speed up the convergence of LSMR, preconditioners with some desirable properties are typically used.
	There are various specific preconditioners, such as incomplete LU \cite{saad:2003}, incomplete QR \cite{lisa:2006},
	and preconditioners based on perturbed QR factorizations \cite{avngto:2009}.
	However, determining whether a given preconditioner is suitable for a particular problem at hand is not straightforward.
	%to determine whether or not a given preconditioner is suitable for a particular
	%problem at hand is not so straightforward.
	%If not effective, many attempts are needed and switched periodically among them.
	
	In 1993, Saad proposed a flexible GMRES (FGMRES) \cite{saad:1993} which still has an Arnoldi-like process, like original GMRES. In the
	process, each matrix-vector multiplications involves a linear system solving that can be viewed as an application
	of some preconditioner that differs from one matrix-vector multiplication to another.
	%multiplication
	%which can be regarded as an inner-outer iterative scheme.
	%Its each outer iterative step involves a linear system
	%setp with each inner capable of  switching preconditioners in the outer iteration.
	%	Specifically, this approach involves solving a  linear system approximately in the inner iteration.
	%	It uses different preconditioners in each outer iteration as a result of different solvers in this inner iteration.
	Accelerating techniques to generate better approximate search space \cite{baje:2005,imli:2016} for GMRES can also be extended to FGMRES, resulting in variants of the method, such as in \cite{Gir:2010,Yang:2022}.  These approaches aim at solving linear systems.
	In 2015, Morikuni and Hayami proposed an inner-outer iterative GMRES method \cite{moha:2015} for solving \eqref{eq:LSP}, where the inner iterations are some stationary iterative methods like NR-SOR and NE-SOR
	to solve normal-equation-type equations in the form of $A^{\T}Av=A^{\T}p$ or $AA^{\T}v=p$.

	In this paper, we will combine the  two ideas above to form a flexible modified LSMR but use non-stationary methods
	to deal with normal-equation-type equations in the inner iteration. Every time the use of a non-stationary method
	yields a preconditioner in 	the Golub-Kahan bidiagonalization process.
	Previously, there are two linear systems to solve
	per iterative step for the right-preconditioned least squares problems, and that can be too demanding computationally.	
	We adopt the concept of factorization-free LSQR (MLSQR) from \cite{arbeha:2014}
	and merge the two linear systems into one to reduce computational cost. We mainly focus on accelerating LSMR rather than LSQR because LSMR exhibts better convergence properties than LSQR.
	However, it's not difficult to apply the same strategy of FMLSMR to create flexible modified LSQR.

	In \cite{chga:2019}, Chuang and Gazzola proposed flexible LSMR (FLSMR) and flexible LSQR for $\ell_p$ regularization
	based on the flexible Golub-Kahan process, in which one upper Hessenberg and one upper triangular matrices are constructed
	via Arnoldi process. Specifically, in FLSMR, two Arnoldi processes are required and so the computational cost is high for
	large scale problems when a really long Arnoldi process is needed due to orthogonalization. However, in our proposed FMLSMR and LSMR,
	only one lower-bidiagonal matrix is generated via the Golub-Kahan process, a two-term recurrence, which
	keeps orthogonalization cost per step low and constant. This is a major difference between FMLSMR and FLSMR.
	%In every step of the Golub-Kahan process, the new orthogonal vectors are generated from  the very previous two vectors. Therefore, our method has much less computational cost in each iteration than that of FLSMR.
	
	The rest of this paper is organized as follows.
	In Section~\ref{sec:2}, we present the flexible modified LSMR method based on modified LSMR for right-preconditioned
	least squares problems and give some theoretical analysis of MLSMR and flexible LSMR.
	In Section~\ref{sec:FMLSMR}, we present the framework of flexible MLSMR with implementation.
	Numerical experiments are shown in Section~\ref{sec:NumTests}. Finally, the conclusion is drawn in Section~\ref{sec:5}.
	
	{\bf{Notation.}} $\bbR^{m\times n}$ is the set of all $m\times n$ real matrices, $\bbR^n=\bbR^{n\times 1}$.
	$I_n$ or $I$ (if its size is clear) is the $n\times n$ identity matrix, and $e_j$ is its $j$th column. $(\cdot)^{\T}$
	takes the transpose of a matrix or vector.
	%Denoted by $i:j$, it's the set of integers from $i$ to $j$ inclusively.
	%	For a matrix $B\in\bbR^{m\times n}$, $$.
	For a vector $u\in\bbR^{n}$, $u_{(i)}$ is its $i$th entry, and $\|\cdot\|_2$ is either $\ell_2$-vector norm or the matrix spectral
	norm:
	$$
	\|u\|_2=\sqrt{\sum\limits_{i} |u_{(i)}|^2},\quad
	\|B\|_2:=\max\limits_{v\neq 0}\frac{\|Bv\|_2}{\|v\|_2}.
	%	\langle u, u\rangle=u^{\T}u=\|u||_2^2,
	$$
	The standard inner product $\langle u,v\rangle=u^{\T}v$ for vectors $u$ and $v$ of the same size, and in particular
	$\langle u, u\rangle=u^{\T}u=\|u||_2^2$.
	Positive definite symmetric matrix $M\in\bbR^{n\times n}$ induces the $M$-inner product
	$\langle u, u\rangle_M=u^{\T}Mu$.
	Finally, the $k$th Krylov subspace $\mathcal{K}_k\left(X,u\right)$ of $X\in \mathbb{R}^{n\times n}$ on $u$ is defined as
	\[
	\mathcal{K}_k\left(X,u\right)=\subspan \{u,Xu,X^2u, \cdots, X^{k-1}u\},
	\]
	spanned by vectors $u,Xu,X^2u, \cdots, $ and $X^{k-1}u$.  $\cR(X)$ is the column space of $X$, spanned by its column vectors.

	\section{Modified LSMR and Flexible LSMR }\label{sec:2}
	In this section, firstly, based on the modified LSQR in \cite{arbeha:2014}, the modified LSMR method (MLSMR) for right-preconditioned least squares problems is introduced and some theoretical analysis of MLSMR is shown. Secondly, we review the flexible LSMR  method (FLSMR) \cite{chga:2019}.
	Lastly, we compare MLSMR and FLSMR from the aspects of theoretical analysis and implementation.
	%\marginpar{\tiny Is MLSMR new?}
	
	\subsection {Modified LSMR}\label{sec:2.1}
	The right-preconditioned least squares problem of \eqref{eq:LSP} is as follows
	\begin{equation}\label{eq:RPLS}
		\min_{\hat{x}}\left\| AL^{-1}\hat{x}-b\right\|_2,\,\, x=L^{-1}\hat{x},
	\end{equation}
	where $L$ is a preconditioner. Any solution to \eqref{eq:RPLS} is a solution to the split-preconditioned normal equation \cite{ake:1996}
	\begin{equation}\label{eq:SPNE}
		L^{-\T}A^{\T}AL^{-1}\hat{x}=L^{-\T}A^{\T}b,
	\end{equation}
	and vice versa.
	
	Recall the Golub-Kahan bidiagonalization process \cite{goka:1965},
	also known as the Lanczos bidiagonalization \cite{ake:1996},
	for a rectangular matrix. It iteratively transforms a matrix into a lower bidiagonal matrix.
	With a right-preconditioner $L$ as in \eqref{eq:RPLS}, the Golub-Kahan bidiagonalization
	process for $AL^{-1}$
	on\footnote {Without loss of generality, we assume initial guess $x_0=0$ in association
		with the least squares problem \eqref{eq:LSP}; otherwise we can always reset $b$ to $b-Ax_0$.}
	$b$,
	is outlined in Algorithm~\ref{alg:GK}.
	\begin{algorithm2e}[t]
		\caption{Golub-Kahan bidiagonalization for $AL^{-1}$}\label{alg:GK}
		\begin{algorithmic}[1]
			\REQUIRE $A\in \bbR ^{m\times n}$, $L\in \bbR ^{n\times n}$, $b\in\bbR^{m}$, and tolerance $\epsilon$; %$1\le k\le\min\{m,n\}$;
			\ENSURE{Partial bidiagonalization of $AL^{-1}\approx U_{k+1}B_kV^T_k$.}
			\vspace{0.5em}
			\STATE {$\beta_1=\left\|b\right\|_2$, $u_1=b/\beta_1$;}
			\STATE {solve $L^{\T}\hat u=A^{\T}u_1$ for $\hat u$;}
			\STATE {$\alpha_1=\|\hat u\|_2$, $v_1=\hat u/\alpha_1$;}
			\FOR{$k=1,2,\ldots,k_{\max}$}
			\STATE{solve $L\tilde v_k=v_k$ for $\tilde{v}_k$;}
			\STATE{ $w=A\tilde v_k$, $q=w-\alpha_ku_k$, $\beta_{k+1}=\|q \|_2$;}
			\IF{$\beta_{k+1}\leq \epsilon\,\|w\|_2$}
			\STATE{{\bf{break}};}
			\ENDIF
			\STATE{$u_{k+1}=q/\beta_{k+1}$;}
			\STATE{solve $L^{\T}\hat u=A^{\T}u_{k+1}$ for $\hat u$;}
			\STATE{$p=\hat u-\beta_{k+1}v_{k}$, $\alpha_{k+1}=\|p \|_2$;}
			\IF {$\alpha_{k+1}\leq \epsilon\,\|\hat u\|_2$}
			\STATE {{\bf {break}};}
			\ENDIF
			\STATE{$v_{k+1}=p/\alpha_{k+1}$;}
			\ENDFOR
			\RETURN{the last $U_{k+1}$, $B_k$, $V_k$.}
		\end{algorithmic}
	\end{algorithm2e}
	%Recall that the Golub-Kahan process with $t=b$ iteratively transforms $[b \quad A]$ to the upper-bidiagonal form $[\beta_1e_1 \quad B_k]$, where
	%
	If Algorithm~\ref{alg:GK} is executed without any breakdown, i.e.,
	all $\alpha_k>0$ and $\beta_k>0$ for $1\leq k \leq k_{\max}$, then
	we will have in theory
	\begin{subequations}\label{eq:PPreGolub}
		\begin{align}
			AL^{-1}V_k&=U_{k+1}B_k, \label{eq:PreGolub_1}\\
			L^{-\T}A^{\T}U_{k+1}& = V_{k}B^{\T}_{k}+\alpha_{k+1} v_{k+1}e^{\T}_{k+1},\label{eq:PreGolub_2}
		\end{align}
	\end{subequations}
	where
	\begin{gather*}
		B_k={\begin{bmatrix}
				\alpha_1 &          &        &        & \\
				\beta_2  & \alpha_2 &        &        &\\
				& \ddots   & \ddots &        & \\
				&          & \beta_k&\alpha_k& \\
				&          &        &\beta_{k+1}&\\
		\end{bmatrix}}, \\
		V_{k}=[v_1, v_2,\cdots, v_k]\in \mathbb{R}^{n\times k},\,\,
		U_{k}=[u_1, u_2,\cdots, u_k] \in \mathbb{R}^{m\times {k}}.
	\end{gather*}
	Both $V_{k}$ and $U_{k}$ are orthonormal.
	Algorithm~\ref{alg:GK} reduces to the original Golub-Kahan process when $L=I$.
	
	It is clear from Algorithm~\ref{alg:GK} that there are two linear systems in the form
	$Lv=w$ and $L^{\T}w=u$ to solve per for-loop.
	Our goal in what follows  is to merge the two linear systems into one in the form
	$Mv=u$, where $M=L^{\T}L$. This does not help when $L$ itself is structured such as
	being lower triangular for which case $Mv=u$ is solved directly via $L^{\T}w=u$ followed by  $Lv=w$,
	but in the case when both $Lv=w$ and $L^{\T}w=u$ have to be solved iteratively, merging
	two linear systems into one per for-loop can bring big savings, not to mention the situation
	when only $M$ is known but not $L$.
	
	Pre-multiplying both sides of (\ref{eq:PreGolub_2}) by $L^{-1}$, we get
	\begin{subequations}\label{eq:AVUB-all}
		\begin{align}
			A\wtd V_{k} &= U_{k+1}B_{k},\label{eq:AVUB}\\
			M^{-1}A^{\T}U_{k+1} &= \wtd V_{k}B^{\T}_{k}+\alpha_{k+1} \tilde{v}_{k+1}e^{\T}_{k+1},\label{eq:MAUVB}
		\end{align}
	\end{subequations}
	where $\wtd V_k=L^{-1}V_k$.
	Notice $\alpha_{k+1}=\|p\|_2=\sqrt{\langle p, p\rangle }$ at Line 12 of Algorithm~\ref{alg:GK}.
	%Here and in what follows, $\langle x, y\rangle:=x^{\T}y$ denotes the standard inner product and
	%$\langle x, y\rangle_M:=x^{\T}M y$ the $M$-inner product.
	It can be verified that
	$$
	\langle p, p\rangle =\langle M^{-1}A^{\T}u_{k+1}-\beta_{k+1}\tilde{v}_{k},
	\,\, M^{-1}A^{\T}u_{k+1}
	-\beta_{k+1}\tilde{v}_{k}\rangle_{M},
	$$
	upon using $M=L^{\T}L$.
	Denote by $s=M^{-1}A^{\T}u_{k+1}-\beta_{k+1}\tilde{v}_{k}$. Then, in Algorithm~\ref{alg:GK},
	$$
	\beta_{k+1}u_{k+1} =A\tilde{v}_{k}-\alpha_{k}u_{k}, \,\,
	\alpha_{k+1} =\sqrt{\langle s, s \rangle_{M}},\,\,
	\tilde{v}_{k+1} =s/\alpha_{k+1}.
	$$
	Define $\tilde p_k=M\tilde{v}_{k}$ and $\tilde{s}=A^{\T}u_{k+1}-\beta_{k+1}\tilde p_k$. We have
	\begin{align*}
		M^{-1}\tilde{s}&=M^{-1}A^{\T}u_{k+1}-\beta_{k+1}M^{-1}\tilde p_k\\
		&=M^{-1}A^{\T}u_{k+1}-\beta_{k+1}\tilde{v}_{k} \\
		&=s.
	\end{align*}
	Hence,
	$\alpha_{k+1}^2=s^{\T}M s=s^{\T}M M^{-1}\tilde{s}=\langle s,\tilde{s}\rangle$.
	Lines 6-16 can be restated as
	\begin{equation*}% \label{PreGolub0}%\left.
		\begin{aligned}
			\beta_{k+1}u_{k+1}&=A\tilde{v}_{k}-\alpha_{k}u_{k},\\
			\tilde{s}&=A^{\T}u_{k+1}-\beta_{k+1}\tilde p_k,\\
			s&=M^{-1}\tilde{s},\\
			\alpha_{k+1}&=\sqrt{\langle s,\tilde{s}\rangle},\\
			\tilde{v}_{k+1}& =s/\alpha_{k+1}.
		\end{aligned}
	\end{equation*}
	It is understood, throughout the paper, that an expression like $s=M^{-1}\tilde{s}$ is really about solving
	linear system $Ms=\tilde{s}$ for $s$.
	%By definition, $\tilde p_{i+1}=M\tilde{v}_{i+1}$.
	Since $\tilde{v}_{k+1} =s/\alpha_{k+1}$ and $s=M^{-1}\tilde{s}$, we have
	\begin{equation*}
		\tilde p_{k+1}=M\tilde{v}_{k+1}=Ms/\alpha_{k+1}=\tilde{s}/\alpha_{k+1}.
	\end{equation*}
	Therefore, a more computationally  cost-effective version of the preconditioned Golub-Kahan bidiagonalization process can be summarized as
	\begin{subequations}\label{PreGolub}
		\begin{align}
			\beta_{k+1}u_{k+1}&=A\tilde{v}_{k}-\alpha_{k}u_{k}, \label{PreGolub-1}\\
			\tilde p&=A^{\T}u_{k+1}-\beta_{k+1}\tilde p_k, \label{PreGolub-2}\\
			\tilde{v}_{k+1}&=M^{-1}\tilde p,  \label{PreGolub-3}\\
			\alpha_{k+1}&=\sqrt{\langle \tilde{v}_{k+1},\tilde p\rangle}, \label{PreGolub-4}\\
			\tilde{v}_{k+1}& =\tilde{v}_{k+1}/\alpha_{k+1}, \label{PreGolub-5}\\
			\tilde p_{k+1} &=\tilde p/\alpha_{k+1}. \label{PreGolub-6}
		\end{align}
	\end{subequations}
	This new formulation has only one linear system
	$M\tilde{v}_{k+1}=\tilde p$ from \eqref{PreGolub-3} to solve.
	
	It can be seen that $U_{k+1}$ is orthonormal, while $\wtd V_{k+1}$ isn't.
	For any $\hat x\in\cR (V_{k})$, $\hat x=V_k\hat y$ for some $\hat{y}\in\bbR^k$.
	Let
	$x=L^{-1}\hat x$ and recall $r=b-Ax$. Using (\ref{eq:AVUB}), we get
	%the residual of the preconditioned normal equation (\ref{eq:SPNE})
	%can be written as
	\begin{align*}
		L^{-\T}A^{\T}r  &=L^{-\T}A^{\T}(b-AL^{-1}\hat{x})=L^{-\T}A^{\T}(b-AL^{-1}V_{k}\hat{y})\\
		&=L^{-\T}A^{\T}(b-U_{k+1}B_k\hat y)=L^{-\T}A^{\T}b-L^{-\T}A^{\T}U_{k+1}B_{k}\hat y\\
		&=\munderbar{\beta}_1v_1-V_{k+1}{\begin{bmatrix}
				B^{\T}_kB_k \\
				\munderbar{\beta}_{k+1}e^{\T}_{k}
			\end{bmatrix}\hat y}\\
		&=V_{k+1}\left(\munderbar{\beta}_1e_1-{\begin{bmatrix}
				B^{\T}_kB_k \\
				\munderbar{\beta}_{k+1}e^{\T}_{k}
			\end{bmatrix}\hat y}\right),
	\end{align*}
	where $\munderbar{\beta}_1=\alpha_1\beta_1$ and $\munderbar{\beta}_{k+1}=\alpha_{k+1}\beta_{k+1}$.
	For the preconditioned least squares problem \eqref{eq:RPLS},
	LSMR seeks an approximate solution $\hat{x}_k\in\cR (V_{k})$
	that minimizes  $\|L^{-\T}A^{\T}r\|_2$ of the preconditioned normal equation (\ref{eq:SPNE})
	over all $x\in\cR (V_{k})$.
	%Recall that i.e., , where $\hat{y}_k$ is the vector we need to find,
	Because $V_{k+1}$ is orthonormal,
	for $\hat x=V_k\hat y$ we have
	\begin{equation}\label{eq:reMLSMR}
		\min_{\hat{x}\in\cR (V_{k})}\left\|L^{-\T}A^{\T}r \right\|_2
		=\min_{\hat{y}}\left\|\munderbar{\beta}_1e_1-{\begin{bmatrix}
				B^{\T}_kB_k \\
				\munderbar{\beta}_{k+1}e^{\T}_{k}
			\end{bmatrix}\hat{y}}\right\|_2,
	\end{equation}
	and hence $\hat{x}_{k}=V_{k}\hat{y}_{k}$, where
	\begin{equation}\label{eq:subMLSMR}
		\hat{y}_{k} =\arg \min\limits_{\hat{y}}\left\|\munderbar{\beta}_1e_1-{\begin{bmatrix}
				B^{\T}_kB_k \\
				\munderbar{\beta}_{k+1}e^{\T}_{k}
			\end{bmatrix}\hat{y}}\right\|_2
	\end{equation}
	which is a least squares problem that can be
	solved by the double QR factorization as in LSMR. Finally,
	%Once the solution of the subproblem is found,
	the solution of the original least squares problem can then be approximated by	
	\begin{equation}\label{eq:subMLSMR'}
		x_{k} = L^{-1} \hat{x}_{k}= L^{-1} V_{k}\hat{y}_{k}=\wtd V_{k}\hat{y}_{k}.
	\end{equation}
	Putting all together, we have established the modified LSMR (MLSMR) method outlined
	in Algorithm~\ref{alg:MLSMR}.
	It should be noted that  Lines 8--19 use the double QR decomposition,
	taken from LSMR \cite{fosa:2011}.

	%\marginpar{\tiny in Alg.~\ref{alg:MLSMR}, how is $\xi_k$ iterated?}

	\begin{algorithm2e}[t]
		\caption{Modified LSMR (MLSMR)}
		\label{alg:MLSMR}
		%\hrule
		\begin{algorithmic}[1]
			\REQUIRE $A\in \bbR^{m\times n}$, $b\in \bbR^{m}$, $M\in\bbR^{n\times n}$, and tolerance $\epsilon$;
			%an integer $k$ with $1\le k\le\min\{m,n\}$;
			\ENSURE  Approximate solution to \eqref{eq:LSP}.
			\vspace{0.5em}
			\STATE {
				$\beta_1=\|b\|_2, u_1=b/\beta_1,\tilde p=A^{\T}u_1,\tilde{v}_1=M^{-1}\tilde p,\alpha_1=\langle \tilde{v}_1,\tilde p \rangle^{1/2}, \tilde p=\tilde p/\alpha_1,\tilde{v}_1=\tilde{v}_1/\alpha_1$,\\
				$\munderbar{\alpha}_{1}=\alpha_1$, $\munderbar{\xi}_{1}=\alpha_{1}\beta_{1}$, $\rho_0=\munderbar{\rho} _{0}=\munderbar{c}_{0}=1$,
				$\munderbar{s}_{0}=0, h_1=\tilde{v}_1$, $\munderbar{h}_{0}=0$, $x_{0}=0$;}
			\FOR{$k=1,2,\ldots,k_{\max}$}
			\STATE  {$\hat u_{k+1}=A\tilde{v}_k-\alpha_ku_k, \beta_{k+1}=\|\hat u_{k+1}\|_2, u_{k+1}=\hat u_{k+1}/\beta_{k+1}$;}
			\STATE  {$\tilde p=A^{\T}u_{k+1}-\beta_{k+1}\tilde p$;}
			\STATE  {$\tilde{v}_{k+1}=M^{-1}\tilde p$, $\alpha_{k+1}=\langle \tilde{v}_{k+1},\tilde p \rangle^{1/2}$;}
			\STATE  {$\tilde p=\tilde p/\alpha_{k+1}$;}
			\STATE {$\tilde{v}_{k+1}=\tilde{v}_{k+1}/\alpha_{k+1}$;}
			\STATE  {$\rho_{k}=(\munderbar{\alpha}^{2}_{k}+\beta^{2}_{k+1})^{1/2}$;}
			\STATE  {$c_{k}=\munderbar{\alpha}_{k}/\rho_{k}$;}
			\STATE  {$s_{k}=\beta_{k+1}/\rho_{k}$, $\theta_{k+1}=s_{k}\alpha_{k}$;}
			\STATE  {$\munderbar{\alpha}_{k+1}=c_{k}\alpha_{k+1}$;}
			\STATE {$\munderbar{\theta}_{k}=\munderbar{s}_{k-1}\rho_{k}$;}
			\STATE  {$\munderbar{\rho}_{k}=((\munderbar{c}_{k-1}\rho_{k})^2+\theta_{k+1}^{2})^{1/2}$;}
			\STATE {$\munderbar{c}_{k}=\munderbar{c}_{k-1}\rho_{k}/\munderbar{\rho}_{k}$;}
			\STATE  {$\munderbar{s}_{k}=\theta_{k+1}/\munderbar{\rho}_{k}$;}
			\STATE  {$\xi_{k}=\munderbar{c}_{k}\munderbar{\xi}_{k}$, $\munderbar{\xi}_{k+1}=-\munderbar{s}_{k}\munderbar{\xi}_{k}$;}
			\STATE {$\munderbar{h}_{k}=h_{k}-(\munderbar{\theta}_{k}\rho_{k}/(\rho_{k-1}\munderbar{\rho}_{k-1}))\munderbar{h}_{k-1}$;}
			\STATE  {$x_{k}=x_{k-1}+(\xi_{k}/(\rho_{k}\munderbar{\rho}_{k}))\munderbar{h}_{k}$;}
			\STATE {$h_{k+1}=\tilde{v}_{k+1}-(\theta_{k+1}/\rho_{k})h_{k}$;}
			\STATE{$r_k=b-Ax_k$;}
			\IF{$\|A^{\T}r_k\|_2\leq \epsilon\, \|A\|_2(\|b\|_2+\|A\|_2\|x_k\|_2)$}
			\STATE {{\bf{break}};}
			\ENDIF
			\ENDFOR
			\RETURN the last $x_k$.
		\end{algorithmic}
	\end{algorithm2e}

	\begin{theorem}\label{theor:2.1}
		At the $k$th step of MLSMR,
		\[x_{k}\in \mathcal{K}_k(M^{-1}A^{\T}A, M^{-1}A^{\T}b)=L^{-1}\mathcal{K}_k(L^{-\T}A^{\T}AL^{-1}, L^{-\T}A^{\T}b).\]
	\end{theorem}
	
	\begin{proof}
		\iffalse
		When LSMR is used on the split-preconditioned normal equation {\rm(\ref{eq:SPNE})},
		we obtain an approximate solution
		\[
		x_{k}\in \cR({\wtd V_k})= \cR (L^{-1}V_{k}).
		\]
		From Algorithm \ref{alg:GK}, $\cR (V_{k})=\mathcal{K}_k(L^{-\T}A^{\T}AL^{-1}, L^{-\T}A^{\T}b).$		
		Thus,
		$$
		\cR({\wtd V_k})=\cR (L^{-1}V_{k})=L^{-1}\mathcal{K}_k(L^{-\T}A^{\T}AL^{-1}, L^{-\T}A^{\T}b).
		$$
		According to the preconditioned Golub-Kahan bidiagonalization process {\rm (\ref{PreGolub})}, we have
		\begin{align*}
			\beta_{1}u_{1}=b, \quad &\alpha_{1}\tilde{v}_{1}=M^{-1}A^{\T}u_{1},\\
			\beta_{k+1}u_{k+1}=A\tilde{v}_{k}-\alpha_{k}u_{k}, \quad &\alpha_{k+1}\tilde{v}_{k+1}=M^{-1}A^{\T}u_{k+1}-\beta_{k+1}\tilde{v}_{k}.
		\end{align*}
		This means
		$x_{k}\in\cR(\wtd V_{k})=\mathcal{K}_k(M^{-1}A^{\T}A, M^{-1}A^{\T}b)$.
		\fi
		By \eqref{eq:subMLSMR'}, we have
		$$
		x_{k}\in \cR(\wtd V_{k})=\cR (L^{-1}V_{k})= L^{-1}\mathcal{K}_k(L^{-\T}A^{\T}AL^{-1}, L^{-\T}A^{\T}b)=\mathcal{K}_k(M^{-1}A^{\T}A, M^{-1}A^{\T}b),
		$$
		as was to be shown.
	\end{proof}

	\begin{theorem}\label{theor:2.2}
		$\wtd V_{k+1}$ is $M$-orthonormal.
	\end{theorem}
	\begin{proof}
		According to the Golub-Kahan bidiagonalization process, we know both $U_{k+1}$ and $V_{k+1}$ are orthonormal, i.e., $U^{\T}_{k+1}U_{k+1}=I_{k+1}$ and $V^{\T}_{k+1}V_{k+1}=I_{k+1}$.
		Since $\wtd V_{k+1}=L^{-1}V_{k+1}$, we get
		\[
		I_{k+1}= V^{\T}_{k+1} V_{k+1}= \wtd V^{\T}_{k+1} L^{\T}L \wtd V_{k+1}= \wtd V^{\T}_{k+1} M \wtd V_{k+1} =\langle \wtd V_{k+1}, \wtd V_{k+1} \rangle_{M},
		\]
		i.e., $\wtd V_{k+1}$ is $M$-orthonormal.
	\end{proof}

	When $A$ has full column rank, \eqref{eq:LSP} has a unique solution.
	Otherwise, there are infinitely many solutions that yield the minimum value of $\|Ax-b\|_2$.
	For the case that $A^{\T}Ax=A^{\T}b$ with singular $A^{\T}A$, it has been proved \cite{fosa:2011}
	that both LSQR and LSMR return the same minimum-norm solution to the least squares problem \eqref{eq:LSP} at convergence.
	We state these conclusions as follows.
	%Therefore, based on Theorem 4.2 and Corollary 4.3 in \cite{fosa:2011}, we have a corollary as following.
	
	\begin{theorem}[{\cite[Corollary 4.3]{fosa:2011}}]\label{coro:2.4}
		At convergence, LSMR returns the minimum-norm solution to \eqref{eq:LSP}.
	\end{theorem}

	%\vspace{0.2em}
	%According to Theorem \ref{coro:2.4} and Corollary \ref{coro:2.5}, we have the similar conclusion for MLSMR.
	\begin{corollary}\label{coro:2.6}
		Suppose that $\hat{x}_*$ is the solution to \eqref{eq:RPLS} obtained via MLSMR at convergence. Then $\hat{x}_*$ is the minimum-norm solution to {\eqref{eq:RPLS}} and $x_*=L^{-1}\hat{x}_*$ is the minimum-M-norm solution to {\rm \eqref{eq:LSP}}.
	\end{corollary}
	
	\begin{proof}
		Since $\hat{x}_*$ is the solution to \eqref{eq:RPLS} obtained via MLSMR, $\hat{x}_*$ is the minimum-norm solution to \eqref{eq:RPLS} according to Theorem \ref{coro:2.4}. Because $\hat{x}_*=Lx_*$, we have $\|\hat{x}_*\|_2=\|x_*\|_M$. This means $x_*$ is the minimum-$M$-norm solution to \eqref{eq:LSP}.
	\end{proof}

	%\marginpar{\tiny This two theorems say that LSQR=LSMR (?)}
	%\marginpar{\tiny At convergence, this two have the same solution. But in numerical, the approximate solutions are different.}
	
	Now we comment on how to develop a modified LSQR in a similar way.
	With the preconditioned Golub-Kahan bidiagonalization process \eqref{PreGolub}, we can obtain the modified LSQR method (MLSQR), i.e., the factorization-free preconditioned LSQR in \cite{arbeha:2014}, by minimizing $\|r\|_2$ over
	$x\in\cR(\wtd V_{k})$. This yields $x_k=\wtd V_{k}\hat{y}_k^{\text{MLSQR}}$, where
	%The reduced problem in MLSQR is
	\begin{equation}\label{eq:subMLSQR}
		\hat{y}_k^{\text{MLSQR}}= \argmin\limits_{\hat{y}}\|\beta_1e_1-B_k\hat{y}\|_2.
	\end{equation}
	Given a nonsingular preconditioner $L$, $\hat{x}=Lx$ is the solution to (\ref{eq:SPNE}), where $x$ is the solution to the normal equation
	\begin{equation}\label{eq:NE}
		A^{\T}Ax=A^{\T}b.
	\end{equation}
	%where we define $x=L^{-1}\hat{x}$. Since the normal equation (\ref{eq:NE}) is consistent, so is (\ref{eq:SPNE}).
	%\argmin\limits_{x\in {\cR{(\wtd V_k)}}}\|r\|_2}

\begin{theorem}[{\cite[Theorem 4.2]{fosa:2011}}]\label{thero:2.3}
	At convergence, LSQR returns the minimum-norm solution to \eqref{eq:LSP}.
\end{theorem}

\begin{corollary}\label{coro:2.5}
	Suppose that $\hat{x}^{\rm{MLSQR}}$ is the solution to \eqref{eq:RPLS} via MLSQR  at convergence. Then $\hat{x}^{\rm{MLSQR}}$ is the minimum-norm solution to \eqref{eq:RPLS} and $x^{\rm{MLSQR}}=L^{-1}\hat{x}^{\rm{MLSQR}}$ is  the minimum-$M$-norm solution to {\rm \eqref{eq:LSP}}.% If $L$ is orthonormal matrix, $x$ is the minimum-norm solution of {\rm \eqref{eq:LSP}}.
\end{corollary}

\begin{proof}
	Since $\hat{x}^{\text{MLSQR}}$ is the solution to \eqref{eq:RPLS} obtained via MLSQR at convergence, $\hat{x}^{\text{MLSQR}}$ is the minimum-norm solution to \eqref{eq:RPLS} according to Theorem \ref{thero:2.3}. Because $\hat{x}^{\text{MLSQR}}=Lx^{\text{MLSQR}}$, we have $\|\hat{x}^{\text{MLSQR}}\|_2=\|x^{\text{MLSQR}}\|_{M}$. This means $x^{\text{MLSQR}}$ is the minimum-$M$-norm solution to \eqref{eq:LSP}.
\end{proof}

\begin{remark}
	{\rm
		MLSMR minimizes the residual of  the normal equation \eqref{eq:SPNE} by solving the subproblem \eqref{eq:subMLSMR}. If $\munderbar{\beta}_{k+1}=0$ for some $k$, then $\alpha_{k+1}=0$ or $\beta_{k+1}=0$. In this case, \eqref{eq:subMLSMR} becomes
		$\min\limits_{\hat{y}}\|\munderbar{\beta}_1e_1-B^{\T}_{k}B_{k}\hat{y}\|_2$ and $B^{\T}_{k}B_{k}\hat{y}_k=\munderbar{\beta}_1e_1$ because  $B_k$ has full rank.
		%It is the normal equation of the least squares problem $\min\limits_{\hat{y}}\|\beta_1e_1-B_k\hat{y}\|$. This means that
		Hence, at  convergence MLSMR returns the same solution as by MLSQR. When $L=I$,  MLSMR reduces to LSMR.
	}
\end{remark}

%\marginpar{\tiny MLSMR=LSMR if $L=I$(?)}

\subsection{Flexible LSMR Method}\label{sec:2.2}
In \cite{chga:2019}, Chung and Gazzola proposed two flexible Krylov subspace methods for $\ell_p$ regularization:
the flexible LSQR (FLSQR) and the flexible LSMR (FLSMR).
The basic idea is to apply FGMRES to construct a flexible variant of the Golub-Kahan process (FGK).
The framework of FGK is shown in Algorithm~\ref{alg:FGK}, where $N_k$ is the preconditioner at the $k$th iteration
determined at runtime.
%\Red{In \cite{chga:2019}, $N_k=?$}

%\marginpar{\tiny $N_k$, $L_k$? You also use $M_k$ later.}
\begin{algorithm2e}[h]
	\caption{Flexible Golub-Kahan Process (FGK)}
	\label{alg:FGK}
	%\hrule
	\begin{algorithmic}[1]
		\REQUIRE $A\in \bbR^{m\times n}$, $b\in \bbR^{m}$,  and tolerance $\epsilon$;
		\ENSURE  $Z_{k}$, $\wtd U_{k+1}$, $W_{k+1}$, $T_{k+1}$, and $P_{k}$.
		\vspace{0.5em}
		\STATE {$\tilde{\beta}_1=\|b\|_2,\tilde{u}_1=b/\tilde{\beta}_1, w_0=0$;}
		\FOR{$k=1, 2, \dots,k_{\max}$}
		\STATE{$v=A^{\T}\tilde{u}_k$, $\xi=\|v\|_2$;}
		\FOR {$j=1, 2, \dots, k-1$}
		\STATE{$t_{j,k}=v^{\T}w_j$;}
		\STATE{$v=v-t_{j,k}w_j;$}
		\ENDFOR	
		\STATE{$t_{k,k}=\|v\|_2$;}
		\IF{$t_{k,k}\leq \epsilon\,\xi $}
		\STATE{\bf break};
		\ENDIF
		\STATE{$w_k=v/t_{k,k};$}
		%\STATE{$v=v-\sum_{j=1}^{k-1}t_{j,k}h_j, t_{k,k}=\|v\|_2,w_k=v/t_{k,k}$;}
		\STATE{solve $N_kz_k=w_k$ for $z_k$ where $N_k\in\bbR^{n\times n}$ is some varying preconditioner determined at runtime;}
		\STATE{$v=Az_k$, $\xi=\|v\|_2$;}
		\FOR{$j=1,\dots,k$}
		\STATE{$p_{j,k}=v^{\T}\tilde{u}_j$;}
		\STATE{$v=v-p_{j,k}\tilde{u}_j;$}
		\ENDFOR
		\STATE{$p_{k+1,k}=\|v\|_2$;}
		\IF {$p_{k+1,k}\leq \epsilon\,\xi$}
		\STATE {\bf break};
		\ENDIF
		\STATE {$\tilde{u}_{k+1}=v/p_{k+1,k}$;}
		\ENDFOR
		\RETURN the last $Z_k, \wtd{U}_{k+1}, W_{k+1}, T_{k+1}, P_{k}.$
	\end{algorithmic}
\end{algorithm2e}

After the $k$th iteration, we have
\begin{equation}\label{eq:AZ}
	AZ_k=\wtd U_{k+1}P_k, \qquad A^{\T}\wtd U_{k+1}=W_{k+1}T_{k+1},
\end{equation}
where
\begin{align*}
	Z_k&=[z_1, z_2, \dots, z_k]=[N^{-1}_1w_1, N^{-1}_2w_2, \dots, N^{-1}_k w_k]\in \mathbb{R}^{n\times k},\\
	\wtd U_{k+1}&=[\tilde{u}_1,\tilde{u}_2,\dots,\tilde{u}_{k+1}]\in\mathbb{R}^{m\times{(k+1)}}, \,\, W_{k+1}=[w_1,w_2,\dots,w_{k+1}]\in\mathbb{R}^{n\times(k+1)},\\
	P_{k}&=[p_{i,j}]\in \mathbb{R}^{(k+1)\times k},\,\, T_{k+1}=[t_{i,j}]\in\mathbb{R}^{(k+1)\times(k+1)}.
\end{align*}
It can be seen that both $\wtd U_{k+1}$ and $W_{k+1}$ are orthonormal, and $P_k$ and $T_{k+1}$
are upper Hessenberg and upper triangular matrices, respectively.
The $k$th approximate solution of FLSMR is given by $x^{\rm FLSMR}_k=Z_ky^{\rm FLSMR}_k$,
where
%Because $W_{k+1}$ is orthonormal, the reduced problem can be simplified as
\begin{equation}\label{eq:ReFLSMR}
	y^{\rm FLSMR}_k=\argmin_y\|\beta_1 t_{1,1}e_1- T_{k+1}P_k y\|_2.
\end{equation}
As to the  residual of  the normal equation (\ref{eq:NE}) at $x^{\rm FLSMR}_k$, we have
\begin{align*}
	A^{\T}r_k& \equiv A^{\T}(b-Ax^{\rm FLSMR}_k) \\
	&=\beta_1 t_{1,1} w_1-A^{\T}AZ_ky^{\rm FLSMR}_k\\
	& =\beta_1 t_{1,1} w_1-W_{k+1}T_{k+1}P_ky^{\rm FLSMR}_k\\
	&=W_{k+1}\left(\beta_1 t_{1,1}e_1- T_{k+1}P_k y^{\rm FLSMR}_k\right).
\end{align*}
%	where
%	%Because $W_{k+1}$ is orthonormal, the reduced problem can be simplified as
%	\begin{equation}\label{eq:ReFLSMR}
	%		y^{\rm FLSMR}_k=\argmin_y\|\beta_1 t_{1,1}e_1- T_{k+1}P_k y\|_2.
	%	\end{equation}
Similarly for FLSQR, $x^{\rm FLSQR}_k=Z_ky^{\rm FLSQR}_k$, where
\begin{align*}
	y^{\rm FLSQR}_k=\argmin_y\|\beta_1 e_1-P_k y\|_2.
\end{align*}
It is shown \cite{chga:2019} that  $x^{\rm FLSQR}_k$ obtained at the $k$th step minimizes the residual norm $\|b-Ax\|_2$
over $x\in \cR({Z_k})$, while $x^{\rm FLSMR}_k$ minimizes $\|A^{\T}(b-Ax)\|_2$ over $x\in\cR({Z}_k)$.
This is the theoretical difference between FLSQR and FLSMR.
We outline FLSMR in Algorithm~\ref{alg:FLSMR}.

\begin{algorithm2e}[h]
	\caption{Flexible LSMR}
	\label{alg:FLSMR}
	%\hrule
	\begin{algorithmic}[1]
		\REQUIRE $A\in \bbR^{m\times n}$, $b\in \bbR^{m}$, $N_k\in\bbR^{n\times n}$ for each $k\geq 1$, and tolerance $\epsilon$;
		\ENSURE  Approximate solution $x_k$ to \eqref{eq:LSP}.
		\vspace{0.5em}
		\STATE {$\tilde{\beta}_1=\|b\|_2,\tilde{u}_1=b/\tilde{\beta}_1, w_0=0$;}
		\FOR{$k=1, 2, \dots,k_{\max}$}
		\STATE{Line 3 - 23 in Algorithm~\ref{alg:FGK};}
		\STATE{Solve for $y_k^{\text{FLSMR}}$ from \eqref{eq:ReFLSMR} }
		\STATE{$x_k=Z_ky_k^{\text{FLSMR}}$, $r_k=b-Ax_k$;}
		\IF{$\|A^{\T}r_k\|_2\leq \epsilon\, \|A\|_2(\|b\|_2+\|A\|_2\|x_k\|_2)$}
		\STATE {{\bf{break}};}
		\ENDIF
		\ENDFOR
		\RETURN the last $x_k$.
	\end{algorithmic}
\end{algorithm2e}

%	\marginpar{\tiny state FLSMR}

\subsection{A Brief Comparison of FLSMR and MLSMR}\label{sec:2.3}
In this subsection, we briefly compare FLSMR with MLSMR.
At appearance, FLSMR possibly employs different preconditioners, i.e., different $N_k$ in Algorithm~\ref{alg:FGK} for each $k$,
while MLSMR uses the same preconditioner, i.e., $M$ in Algorithm~\ref{alg:MLSMR}.
When all $N_k$ is taken to be the same as $M$, we have the following result.

\begin{theorem}\label{theor:2.9}
	If the preconditioners $N_k$ in Algorithm~\ref{alg:FGK} are the same as the preconditioner $M=L^{\T}L$ in Algorithm~\ref{alg:MLSMR},
	then the solutions by FLSMR and MLSMR satisfy
	\begin{equation}\label{eq:FLSMR}
		x^{\rm FLSMR}_k=\argmin_x\|A^{\T}(b-Ax)\|_2 \quad \text{over}\quad x\in\mathcal{K}_k(M^{-1}A^{\T}A, M^{-1}A^{\T}b),
	\end{equation}
	\begin{equation}\label{eq:MLSMR}
		x^{\rm MLSMR}_k=\argmin_x\|A^{\T}(b-Ax)\|_{M^{-1}}\quad \text{over} \quad x\in\mathcal{K}_k(M^{-1}A^{\T}A, M^{-1}A^{\T}b),
	\end{equation}
	respectively.
\end{theorem}

\begin{proof}
	According to \eqref{eq:AZ}, we have
	\[
	A^{\T}AM^{-1}W_{k}=W_{k+1}T_{k+1}P_{k},
	\]
	where $W_{k}$ is orthonormal  and $T_{k+1}P_{k}$ is upper Hessenberg. Hence,
	FLSMR is exactly the same as GMRES applied to the right-preconditioned linear system
	\[
	A^{\T}AM^{-1}\hat{x}=A^{\T}b, \quad x=M^{-1}\hat{x}.
	\]
	Therefore, at the $k$th iteration, we have
	$$x^{\rm FLSMR}_k=M^{-1}\hat{x}_k\in  M^{-1}\mathcal{K}_k(A^{\T}AM^{-1}, A^{\T}b)=\mathcal{K}_k(M^{-1}A^{\T}A, M^{-1}A^{\T}b).$$
	
	According to Theorem~\ref{theor:2.1}, both FLSMR and MLSMR search their $k$th approximate solutions over
	$$
	x\in\cR{(Z_k)}=\mathcal{K}_k(M^{-1}A^{\T}A, M^{-1}A^{\T}b).
	$$
	This implies \eqref{eq:FLSMR} because $x^{\rm FLSMR}_k$ minimizes $\|A^{\T}(b-Ax)\|_2$ over $\cR({Z}_k)$. On othe other hand by \eqref{eq:reMLSMR}, we know that $x^{\rm MLSMR}_k=L^{-1}\hat{x}_k$, where
	$$
	\hat{x}_k=\argmin_{\hat{x}\in \mathcal{K}_k(L^{-\T}A^{\T}AL^{-1}, L^{-\T}A^{\T}b)}\|L^{-\T}A^{\T}(b-AL^{-1}\hat{x})\|_2.
	$$
	Since $\|L^{-\T}A^{\T}(b-AL^{-1}\hat{x})\|_2^2=\|A^{\T}(b-Ax)\|_{M^{-1}}^2$ and
	$$
	L^{-1}\mathcal{K}_k(L^{-\T}A^{\T}AL^{-1}, L^{-\T}A^{\T}b)=\mathcal{K}_k(M^{-1}A^{\T}A, M^{-1}A^{\T}b),
	$$
	we conclude
	$$
	x^{\rm MLSMR}_k=\argmin_{x\in\cR{(Z_k)}}\|A^{\T}(b-Ax)\|_{M^{-1}},
	$$	
	as was to be shown.
\end{proof}

Theorem~\ref{theor:2.9} states that FLSMR and MLSMR solve different optimization problems over the same search space.

In MLSMR, the $k$th approximate solution $x_k$ satisfies
\begin{equation}\label{eq:ReMLSMR}
	x_k=L^{-1}\hat{x}_k,  \quad \hat{x}_k=
	\argmin_{\hat{x}\in \mathcal{K}_k(L^{-\T}A^{\T}AL^{-1}, L^{-\T}A^{\T}b)}
	\|L^{-\T}A^{\T}r\|_2,
\end{equation}
where $r=b-AL^{-1}\hat{x}$. This involves $L$. In what follows, we will transform it into one involving
$M=L^{\T}L$ only. For any $\hat x\in\cR (V_{k})=\mathcal{K}_k(L^{-\T}A^{\T}AL^{-1})$, $\hat x=V_k\hat y$ for some $\hat{y}\in\bbR^k$, where $V_k$ is orthonormal and obtained in Algorithm \ref{alg:GK}. Based on  the previous discussion about MLSMR in Subsection~\ref{sec:2.1}, we know that
\begin{equation}\label{eq:ReMLSMR.1}
	L^{-\T}A^{\T}r
	=V_{k+1}\left(\alpha_1\beta_1e_1-{\begin{bmatrix}
			B^{\T}_kB_k \\
			\alpha_{k+1}\beta_{k+1}e^{\T}_{k}
		\end{bmatrix}\hat{y}}\right).
\end{equation}
Denote by $D_{k+1}={\begin{bmatrix}
		B^{\T}_kB_k \\
		\alpha_{k+1}\beta_{k+1}e^{\T}_{k}
\end{bmatrix}}$. Multiplying both sides of (\ref{eq:ReMLSMR.1}) by $L^{-1}$, we get
$$
M^{-1}A^{\T}(b-AL^{-1}\hat{x})
=\wtd V_{k+1}\left(\alpha_1\beta_1e_1-D_{k+1}\hat{y}\right),
$$
where $\wtd V_{k+1}=L^{-1}V_{k+1}$, which is defined in Subsection~\ref{sec:2.1}.
Thus,
\begin{equation}\label{eq:ReMLSMR.2}
	A^{\T}r=M\wtd V_{k+1}\left(\alpha_1\beta_1e_1-D_{k+1}\hat{y}\right).
\end{equation}
Let $M\wtd V_{k+1}=\check{Q}_{k+1}\check{R}_{k+1}$ be the QR decomposition of $M\wtd V_{k+1}$,
where $\check{Q}_{k+1}\in\mathbb{R}^{n\times{(k+1)}}$ is orthonormal and $\check{R}_{k+1}\in \mathbb{R}^{{(k+1)}\times{(k+1)}}$ is upper triangular. We then rewrite (\ref{eq:ReMLSMR.2}) as
$$
A^{\T}r=\check{Q}_{k+1}\check{R}_{k+1}\left(\alpha_1\beta_1e_1-D_{k+1}\hat{y}\right).
$$
Naturally, this leads to a new way for the original least squares problem. Namely, instead of (\ref{eq:ReMLSMR}),
we  seek an approximation as follows:
\begin{equation}\label{eq:RePMLSMR}
	\tilde x_k=\argmin_{x\in\mathcal{K}_k(M^{-1}A^{\T}A, M^{-1}A^{\T}b)}\|A^{\T}r\|_2=\wtd V_{k}\hat{z}_k
\end{equation}
%	\marginpar{\tiny $\tilde x_k=V_kz_k$?}
where
$$
\hat{z}_k=\argmin_{\hat{y}}\left\|\alpha_1\beta_1\check{R}_{k+1}e_1-\check{R}_{k+1}D_{k+1}\hat{y}\right\|_2.
$$
Notice that $\alpha_1\check{R}_{k+1}e_1=\alpha_1r_{1,1}e_1=\|A^{\T}u_1\|_2e_1=t_{1,1}e_1$, and
$\check{R}_{k+1}D_{k+1}$ is upper Hessenberg, we get
\begin{equation}\label{eq:ReNMLSMR}
	\hat{z}_k=\argmin_{\hat{y}}\left\|\beta_1t_{1,1}e_1-\check{R}_{k+1}D_{k+1}\hat{y}\right\|_2,
\end{equation}
which takes the same form as FLSMR's reduced problem (\ref{eq:ReFLSMR}).
With the above analysis, we can get a variant of MLSMR, which is
the preconditioned Golub-Kahan process \eqref{eq:PPreGolub} followed by solving \eqref{eq:ReNMLSMR}. We can see that this new variant of MLSMR is equivalent to FLSMR since they both minimize the same objective over the same search space, $x\in\cR(\wtd V_{k})=\cR(Z_k)=\mathcal{K}_k(M^{-1}A^{\T}A, M^{-1}A^{\T}b)$.

\section{Flexible Modified LSMR}\label{sec:FMLSMR}
\subsection{FMLSMR}\label{ssec:FMLSMR}
In MLSMR (Algorithm~\ref{alg:MLSMR}), the most extreme but impractical preconditioner is $M=A^{\T}A$,
% and $M^{-1}$ can be exactly calculated without rounding any error,
with which $x_1$ is the exact solution to the normal equation (\ref{eq:normal_type}).
However, that is not feasible in practice for large $n$. Some approximate inverse of  $A^{\T}A$ has to be used, or, equivalently,
to solve approximately
\begin{equation}\label{eq:normal_type}
	A^{\T}A\tilde{v}=\tilde p,
\end{equation}
for $\tilde{v}$ at Lines 1 and 5 in Algorithm ~\ref{alg:MLSMR}. By doing so, we  implicitly determine some approximations, likely unknown but exist, of $(A^{\T}A)^{-1}$
in the inner iterations.
Using stationary methods such as the Jacobi and SOR-type methods to solve (\ref{eq:normal_type}) yields  preconditioners $M$
that remains the same for each inner iteration, but applying non-stationary methods like CG or MINRES dynamically
selects varying preconditioners,
i.e., $M$ changes as $\tilde p$ changes from one iteration to the next. 
%In fact, CG or MINRES on \eqref{eq:normal_type} with initial $0$ yields
%an approximate solution $\tilde{v}\approx f(A^TA)\tilde p$ with $f(\cdot)$ is some low degree polynomial dependent of $\tilde p$ such that $f(0)=1$, or in the other words $\tilde{v}\approx M^{-1}\tilde p$ for some matrix $M$ dependent of $\tilde p$.
The latter leads to our  new approach, namely the {\em flexible modified LSMR method\/} (FMLSMR),
as outlined in Algorithm~\ref{alg:FMLSMR}.
In our numerical tests later in Section~\ref{sec:NumTests},  MINRES is used to solve \eqref{eq:normal_type}
for FMLSMR.

\begin{algorithm2e}[h]
	\caption{Flexible MLSMR}
	\label{alg:FMLSMR}
	%\hrule
	\begin{algorithmic}[1]
		\REQUIRE $A\in \bbR^{m\times n}$, $b\in \bbR^{m}$, implicitly varying preconditioner $M_k\in\bbR^{n\times n}$,
		and \\tolerance $\epsilon$;
		\ENSURE  Approximate solution to \eqref{eq:LSP}.
		\vspace{0.5em}
		\STATE $\beta_1=\|b\|_2, u_1=b/\beta_1,\tilde p_1=A^{\T}u_1$;
		\STATE solve $A^TA\tilde{v}_1=\tilde p_1$ approximately for $\tilde{v}_1$;
		\STATE $\alpha_1=\langle \tilde{v}_1,\tilde p_1 \rangle^{1/2}$, 
		$\hat p_1=\tilde p_1/\alpha_1,\tilde{v}_1=\tilde{v}_1/\alpha_1$;
		\STATE $\munderbar{\alpha}_{1}=\alpha_1$, $\munderbar{\xi}_{1}=\alpha_{1}\beta_{1}$,
		$\rho_0=\munderbar{\rho}_{0}=\munderbar{c}_{0}=1$,
		$\munderbar{s}_{0}=0, h_1=\tilde{v}_1$, $\munderbar{h}_{0}=0$, $x_{0}=0$;
		\FOR{$k=1,2,\ldots, k_{\max}$}
		\STATE  {$\hat u_{k+1}=A\tilde{v}_k-\alpha_ku_k, \beta_{k+1}=\|\hat u_{k+1}\|_2, u_{k+1}=\hat u_{k+1}/\beta_{k+1}$;}
		\STATE  $\tilde p_{k+1}=A^{\T}u_{k+1}-\beta_{k+1}\hat p_k$;
		\STATE  solve $A^TA\tilde{v}_{k+1}=\tilde p_{k+1}$ approximately for $\tilde{v}_{k+1}$;
		\STATE  $\alpha_{k+1}=\langle \tilde{v}_{k+1},\tilde p_{k+1} \rangle^{1/2}$;
		\STATE  $\hat p_{k+1}=\tilde p_{k+1}/\alpha_{k+1}$, $\tilde{v}_{k+1}=\tilde{v}_{k+1}/\alpha_{k+1}$;
		\STATE  Lines 8-23 of Algorithm \ref{alg:MLSMR};
		\ENDFOR
		\RETURN the last $x_k$.
	\end{algorithmic}
\end{algorithm2e}

%\marginpar{\tiny Alg.~\ref{alg:FMLSMR} won't run with \Red{Lines 10-23} of Algorithm \ref{alg:MLSMR}}

Symbolically, we may write $\tilde{v}_1=M_1^{-1}\tilde p_1$ at Line 2 and $\tilde{v}_{k+1}=M_{k+1}^{-1}\tilde p_{k+1}$
at Line 8, where $M_1$ and $M_{k+1}$ are dependent of vectors $\tilde p_1$ and $\tilde p_{k+1}$, respectively, and of computed approximations $\tilde{v}_1$ and $\tilde{v}_{k+1}$, respectively, as well. Exactly, what these $M_k$ are is not important
as far as executing Algorithm~\ref{alg:FMLSMR} is concerned.
%, but for posterior analysis such as
%in what follows for bounding $\left \|A^{\T}r_k\right \|_2$, we may take any $M_k$ that maps $\tilde v_k$ to $\tilde p_k$, e.g.,
%$$
%M_k=I
%$$
Evidentally, the approximate solution $x_k$ is sought in $\cR(\wtd V_{k})$.
%, which is not the search space $\mathcal{K}_{k}(M^{-1}A^{\T}A, M^{-1}A^{\T}b)$ we mentioned before.
After the $k$th step, we have
\begin{align*}
	A\wtd V_{k}&=U_{k+1}B_{k}, \\
	A^{\T}U_{k+1}
	&= \left[M_{1}\tilde{v}_{1}, M_{2}\tilde{v}_{2}, \cdots, M_{k+1}\tilde{v}_{k+1}\right]
	\begin{bmatrix}
		B^{\T}_k\\
		\alpha_{k+1}\beta_{k+1}e^{\T}_{k+1}
	\end{bmatrix}.
	%&=\left[M_{1}\tilde{v}_{1}, M_{2}\tilde{v}_{2}, \cdots,   M_{k+1}\tilde{v}_{k+1}\right]L^{\T}_{k+1}.\\
\end{align*}
Let $x_k=\wtd V_{k}\tilde{y}_k$. $\|A^{\T}r_k\|_2$ can be expressed as follows:
\begin{align*}
	\left \|A^{\T}r_k\right \|_2
	& =\left \|A^{\T}b-A^{\T}Ax_k \right \|_2=\left \|A^{\T}b-A^{\T}U_{k+1}B_{k}\tilde{y}_k\right \|_2\\
	& =\left \|A^{\T}b-\left[M_{1}\tilde{v}_{1}, \cdots, M_{k+1}\tilde{v}_{k+1}\right]\begin{bmatrix}
		B^{\T}_{k}B_{k}\\                 \alpha_{k+1}\beta_{k+1}e^{\T}_{k} \end{bmatrix}
	\tilde{y}_k\right \|_2\\
	& =\left\|\left[M_{1}\tilde{v}_{1}, \cdots, M_{k+1}\tilde{v}_{k+1}\right](\alpha_1\beta_{1}e_{1}-\begin{bmatrix}
		B^{\T}_{k}B_{k}\\
		\alpha_{k+1}\beta_{k+1}e^{\T}_{k}
	\end{bmatrix}
	\tilde{y}_k)\right\|_2\\
	&\leq\|[\hat{p}_{1}, \cdots, \hat{p}_{k+1}]\|_2
	\left\|\alpha_1\beta_{1}e_{1}-\begin{bmatrix}
		B^{\T}_{k}B_{k}\\
		\alpha_{k+1}\beta_{k+1}e^{\T}_{k}
	\end{bmatrix}\tilde{y}_k\right\|_2,
\end{align*}
where
\[
\tilde{y}_k=\argmin\limits_{\tilde{y}}\left\|\alpha_1\beta_{1}e_{1}-\begin{bmatrix}
	B^{\T}_{k}B_{k}\\
	\alpha_{k+1}\beta_{k+1}e^{\T}_{k}
\end{bmatrix}\tilde{y}\right\|_2.
\]
Hence, if $\|[\hat{p}_{1}, \cdots, \hat{p}_{k+1}]\|_2$ is not too large, $x_{k} $ can be a good approximate solution to the original least squares problem \eqref{eq:LSP}.
Because $[\hat{p}_{1}, \cdots, \hat{p}_{k+1}]$ is not orthonormal,  the residuals by FMLSMR may not be monotonically decreasing, unlike in FLSMR.
%\marginpar{\tiny why?}
However, the computation cost of FMLSMR is less than that of FLSMR because in FLSMR, a Gram-Schmit orthogonalization step is conducted twice at Lines 4 -- 7 and Lines 15 -- 18 of Algorithm~\ref{alg:FGK}.

\section{Numerical Experiments}\label{sec:NumTests}
In this section, we perform numerical tests to demonstrate the advantage of our method FMSLMR.

Firstly, we compare the computational cost of  LSMR, FLSMR, and FMSLMR. Table~\ref{tbl:flops4} lists the numbers
of flops of all methods where the $\ell$-step Lanczos process is used in the inner solver for both FLSMR and FMLSMR.
The symbol  ``\texttt{MV}"  denotes the number of flops required for  a single matrix-vector multiplication with
$A\in \mathbb {R}^{m\times n}$, which is taken to be twice the number of nonzero entries in $A$.
Only the dominant terms are included in flops in each iteration, which are matrix-vector multiplications, solutions
of the inner linear systems and vector-vector operations. Both the computational costs of FLSMR
and FMLSMR are more than LSMR because of their inner iterations, and at the same time, for the same inner solver, the computational cost of
FMLSMR is less than that of FLSMR.
Later we will also report CPU times for all examples.

\renewcommand{\arraystretch}{1.2}
\begin{table}[h]
	\centering
	\caption{Flops of LSMR Variants}
	\label{tbl:flops4}
	\begin{tabular}{|l|l|}
		\hline
		% after \\: \hline or \cline{col1-col2} \cline{col3-col4} ...
		$k$-step of LSMR   & $2k$ (\texttt{MV})  +$(8n+2m)k$\\
		%per cycle of MLSMR(GMRES(k)) ($k$) & $(k+3)$(MV)+$2k^2n+4k^2$  \\
		\hline
		$k$-step of FLSMR($\ell$) & $(2\ell+4)k$ (\texttt{MV})+$8k\ell n+4k\ell^2+2/3k^3+2(m+n+1)k^2$ \\
		\hline
		$k$-step of FMLSMR ($\ell$) & $(2\ell+4)k$ (\texttt{MV})+$8k\ell n+4k\ell^2+(12n+2m)k$  \\
		\hline
	\end{tabular}
\end{table}

Secondly, we compare the storage requirements for the three methods after the $k$th iteration in Table~\ref{tbl:storage}. In the table, the 2nd and 3rd columns show the numbers of stored vectors at the $k$th iteration with respect to different dimensions. The last column displays the stored matrices for each methods. FLSMR is the only one requiring matrices storage.
The storage of FLSMR increases quadratically in $k$. However, FMLSMR consumes the comparable amount of storage as LSMR,
which is far less than that of FLSMR. The results of our numerical examples will confirm this advantage of FMLSMR over FLSMR.

%\marginpar{\tiny variables in Table~\ref{tbl:storage} make sense to the reader?}

\begin{table}[h]
	\caption{Storage of LSMR Variants}
	\label{tbl:storage}
	\centering
	\begin{tabular}{|l|l|l|l|}
		\hline
		% after \\: \hline or \cline{col1-col2} \cline{col3-col4} ...
		&  \# of vectors in $\mathbb{R}^m$  & \# of vectors in $\mathbb{R}^n$  & Matrices \\
		\hline
		$k$-step of LSMR   &  3 & 5 & \\
		\hline
		$k$-step of FLSMR($\ell$) & $k+4$ & $
		k+3$&  $T_{k+1}$, $P_{k+1}$ \\
		\hline
		$k$-step of FMLSMR($\ell$)& 3& 7& \\
		\hline
	\end{tabular}
\end{table}

Third, we report our numerical results on 8 testing problems which are drawn
from the SuiteSparse Matrix Collection\footnote{https://sparse.tamu.edu/}
and Matrix Market\footnote{https://math.nist.gov/MatrixMarket}.
%Before that, 	we explain some details of our examples and how we conduct our numerical tests.
Among the problems, {\tt{crack}}, {\tt{biplane-9}} and {\tt{delaunay\_n16}} come with a right-hand side $b$.
%\marginpar{\tiny {\tt randn} or {\tt rand}?}
A random vector $b$ is generated by {\tt{rand}} for each of the rest of the problems.
Table~\ref{tab:matrices4} lists some of their important characteristics, including the matrix size $m$ and $n$,
the number {\tt{nnz}} of nonzero entries in $A$, and the sparsity {\tt{nnz}}/$(mn)$.
These are  representatives of many other problems from the collections we have tested.
All  tests are done by MATLAB (version R2020b) on a Mac PC with 2.7 GHz Intel Core i7 and 16GB memory.

\begin{table}[ht]
	\caption{Testing Matrices}
	\label{tab:matrices4}
	\centering
	\begin{tabular}{c c c c c c c}
		\hline
		% after \\: \hline or \cline{col1-col2} \cline{col3-col4} ...
		ID & matrix & $m$   &$n$& \texttt{nnz} & sparsity  \\
		\hline
		\hline
		%\rowcolor[gray]{.93}
		1 & \texttt{well1850}     & 1850   & 712  & 8755    & $1.45\times 10^{-2}$     \\
		2 & \texttt{cat\_ears\_3\_4}     & 13271   & 5226  & 39592   & $5.7087\times 10^{-4}$     \\
		%\rowcolor[gray]{.93}
		3 & \texttt{delaunay\_n16}   & 65536  &  65536 & 393150  & $9.1537\times 10^{-5}$     \\
		4 & \texttt{biplane-9}     & 21701   & 21701  & 84076   & $1.7853\times 10^{-4}$     \\
		5 & \texttt{flower\_7\_4}     & 67593  & 27693  & 202218 & $1.0803\times 10^{-4}$     \\
		%\rowcolor[gray]{.93}
		6 & \texttt{crack}      & 10240   &10240 & 60760    & $5.7936\times 10^{-4} $    \\
		7 & \texttt{fe\_body}    & 45087 & 45087 &327468  & $1.6109\times 10^{-4} $   \\
		%\rowcolor[gray]{.93}
		8 & \texttt{stufe-10}    &   24010  & 24010   & 92828   & $1.6103\times 10^{-4}$    \\
		\hline
	\end{tabular}
\end{table}

The stopping criteria are either
\begin{equation}\label{eq:NRes}
	\text{NRes}=\frac{\|A^{\T}(Ax-b)\|_2}{\|A\|_1(\|A\|_1\|x\|_2+\|b\|_2)}\le 10^{-12},
\end{equation}
on normalized residual (NRes) or the number of iterations reaches $10^5$,
where using matrix $\ell_1$-norm $\|A\|_1$ is for its easiness in computation.
Another assessment is the backward error for an approximate least squares solution $x$, which measures the perturbation to $A$ that would make $x$ an exact least squares solution to a perturbed least squares problem:
\begin{equation}\label{eq:BEopt}
	\mu (x)\equiv \min_{E}\|E\|_F \quad {\text{s.t.}}\quad (A+E)^{\T}(A+E)x=(A+E)^{\T}b,
\end{equation}
Wald{\'e}n \cite{Wald:1995}  et al.  and Higham \cite{high:2002} proved that the backward error is the smallest singular value of  the matrix
\begin{equation*}
	\begin{bmatrix}
		A & \frac{\|r\|_2}{\|x\|_2}\left(I-\frac{rr^{\T}}{\|r\|_2^2}\right)
	\end{bmatrix}.
\end{equation*}
In our numerical experiments, we use the following easily computable estimate of backward error in Stewart \cite{stew:1977},
\begin{equation}\label{eq:BErr}
	\hat{E}=-\frac{rr^{\T}A}{\|r\|_2^2}, \quad \|\hat{E}\|_2=\frac{\|A^{\T}r\|_2}{\|r\|_2},
\end{equation}
which satisfies the constraint in \eqref{eq:BEopt} but does not achieve the minimum there.
%The normalized backward error $\|\hat{E}\|_2/\|A\|_1$ is displayed in our results.
Backward error is widely used to estimate the accuracy and stability of a method for least squares problems.
It is usually accepted that the smaller the backward error is, the more accurate the approximate solution is.

\begin{table}[h]
	\caption{Number of Iterations}\label{tbl:cycles4}
	\centering
	\begin{tabular}{|c|c|c|c|c|c|}
		\hline
		ID &{matrix}                  &$\ell$ &LSMR      & FLSMR      &FMLSMR    \\ \hline\hline
		1  &{\tt well1850}         &8  & 463       & 167    & {\bf 117} \\
		2 &{\tt cat\_ears\_3\_4}  & 8 & 163 & 26 & {\bf 25}\\
		3  &\texttt{delaunay\_n16}   &30   & --   &  --  & {\bf 14571} \\
		4  &\texttt{biplane-9}   &15   & -- & --     & {\bf 24611}\\
		5  &\texttt{flower\_7\_4}  &8   & 195   & 29   &{\bf  28}\\
		%				5 & \texttt{cs4}  &30  & 42528    & 1591   &2492\\
		%				6 & \texttt{cti} & 30 & 44621& 1665& 2542\\
		6 & \texttt{crack} &10 &  54015    &   --     &  {\bf{7400}}    \\
		7 & \texttt{fe\_body}  &30 & -- & --   &    {\bf 11200}   \\
		8 & \texttt{stufe-10}  & 10 &  28047&  --& {\bf 3297}\\					
		\hline
	\end{tabular}
\end{table}

\begin{figure}
	{\centering
		\begin{tabular}{cc}
			\hspace{-0.3 cm}
			\resizebox*{0.44\textwidth}{0.235\textheight}{\includegraphics{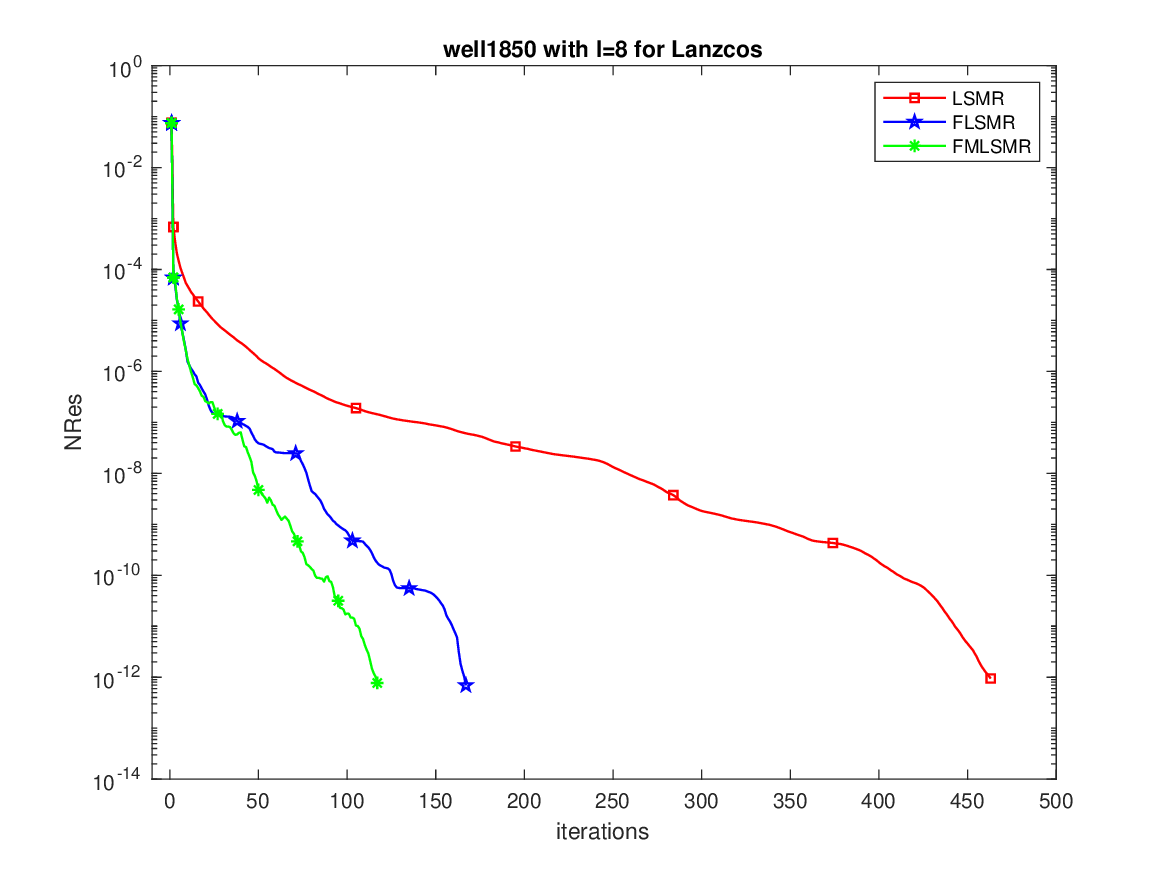}}
			&  \hspace{-0.5 cm}
			\resizebox*{0.44\textwidth}{0.235\textheight}{\includegraphics{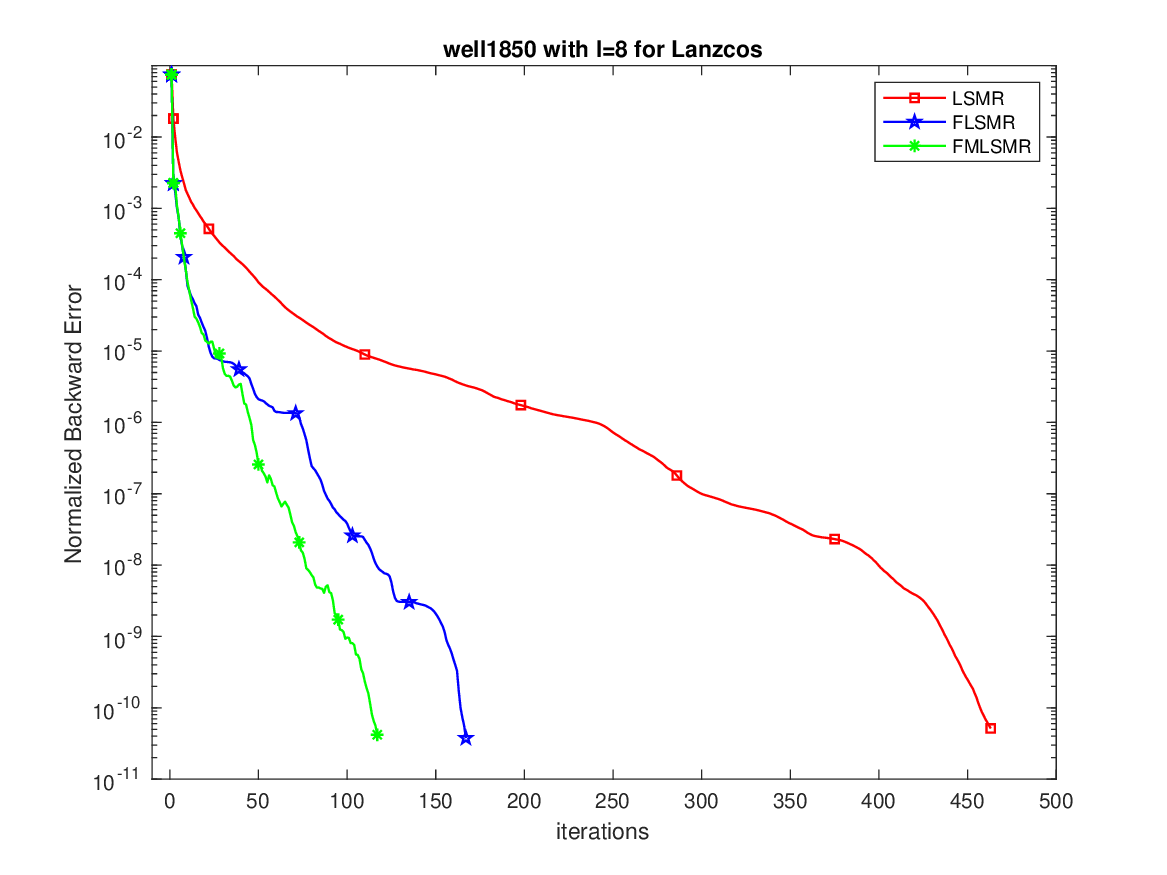}}\\
			\hspace{-0.3 cm}
			\resizebox*{0.44\textwidth}{0.235\textheight}{\includegraphics{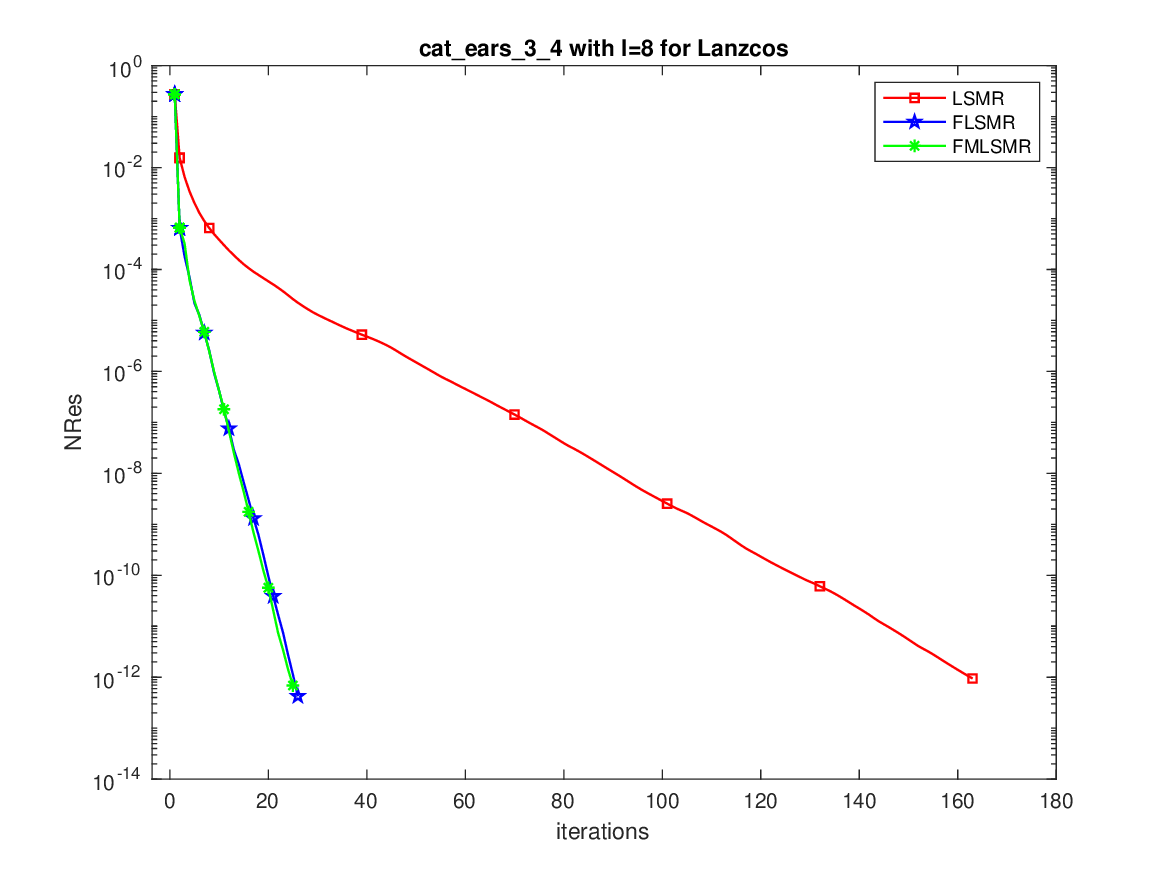}}
			&  \hspace{-0.5 cm}
			\resizebox*{0.44\textwidth}{0.235\textheight}{\includegraphics{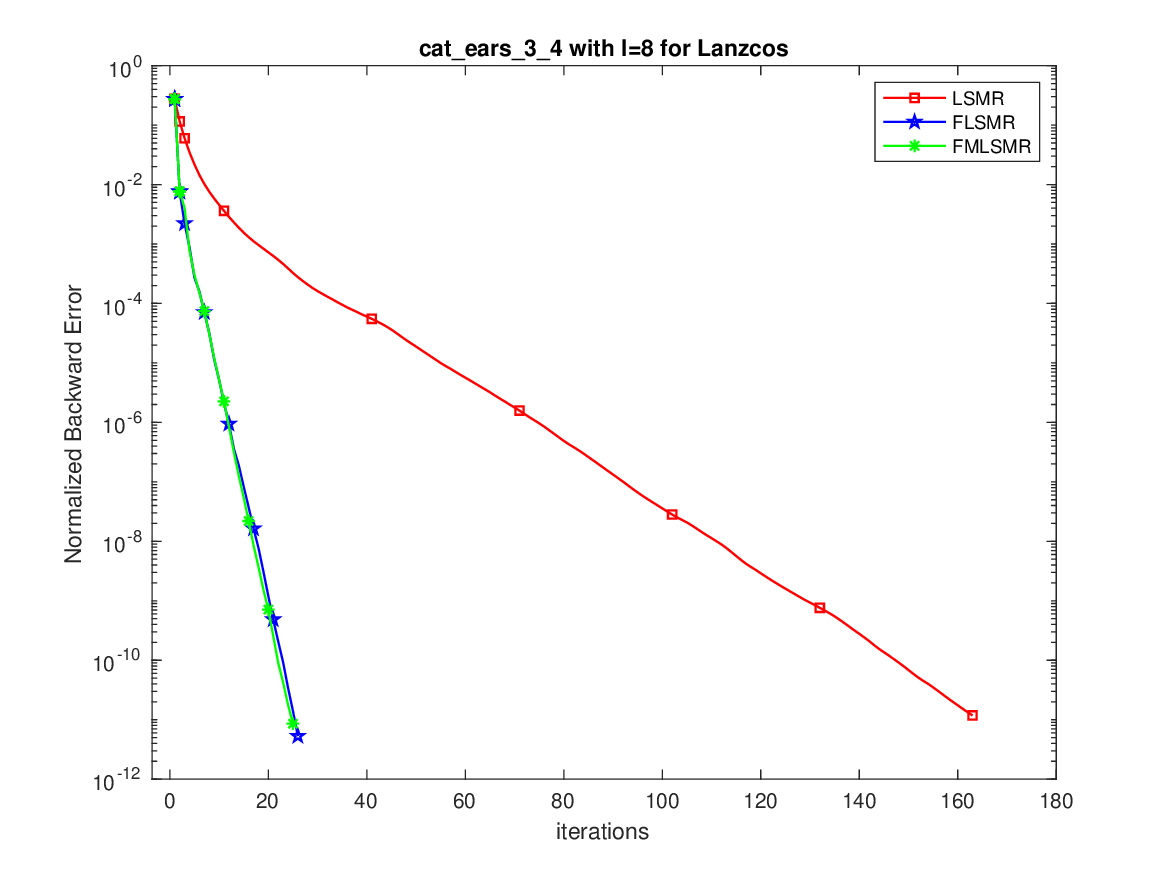}}\\
			\hspace{-0.3 cm}
			\resizebox*{0.44\textwidth}{0.235\textheight}{\includegraphics{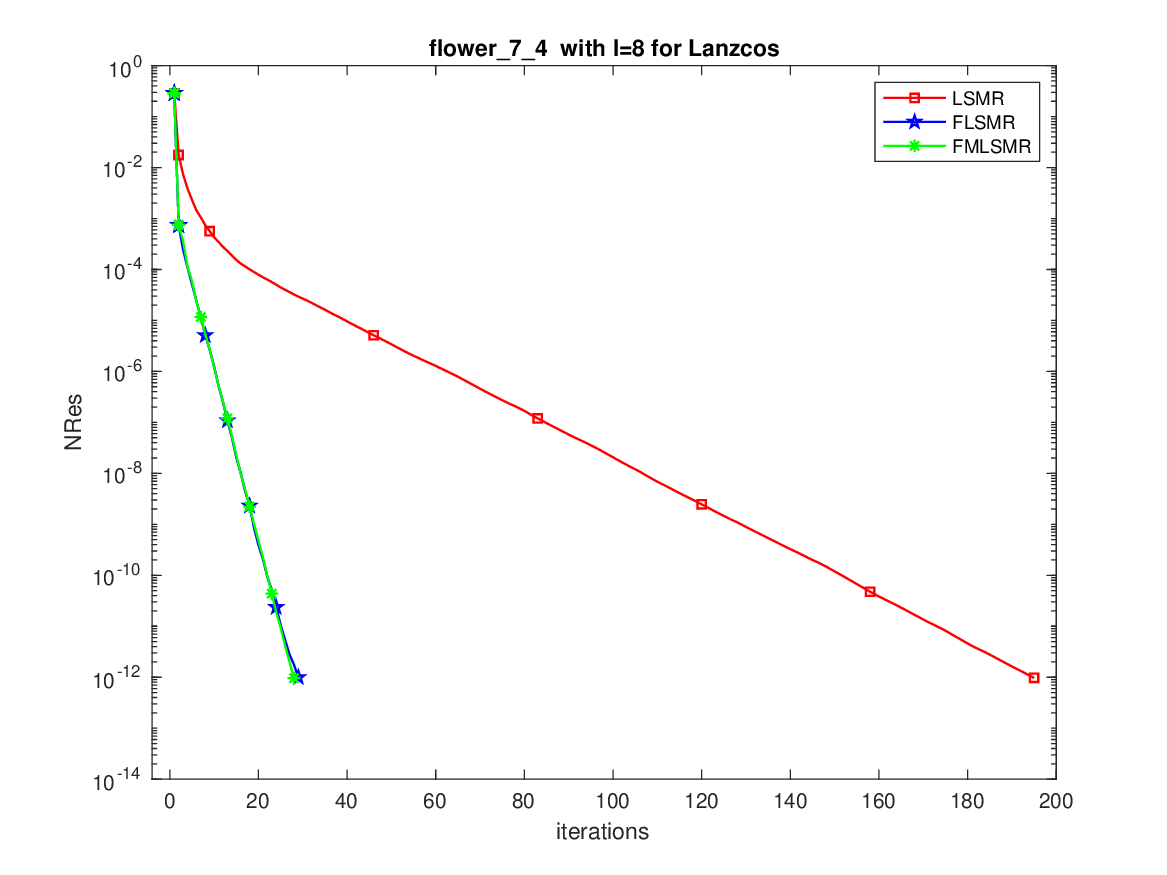}}
			&  \hspace{-0.5 cm}
			\resizebox*{0.44\textwidth}{0.235\textheight}{\includegraphics{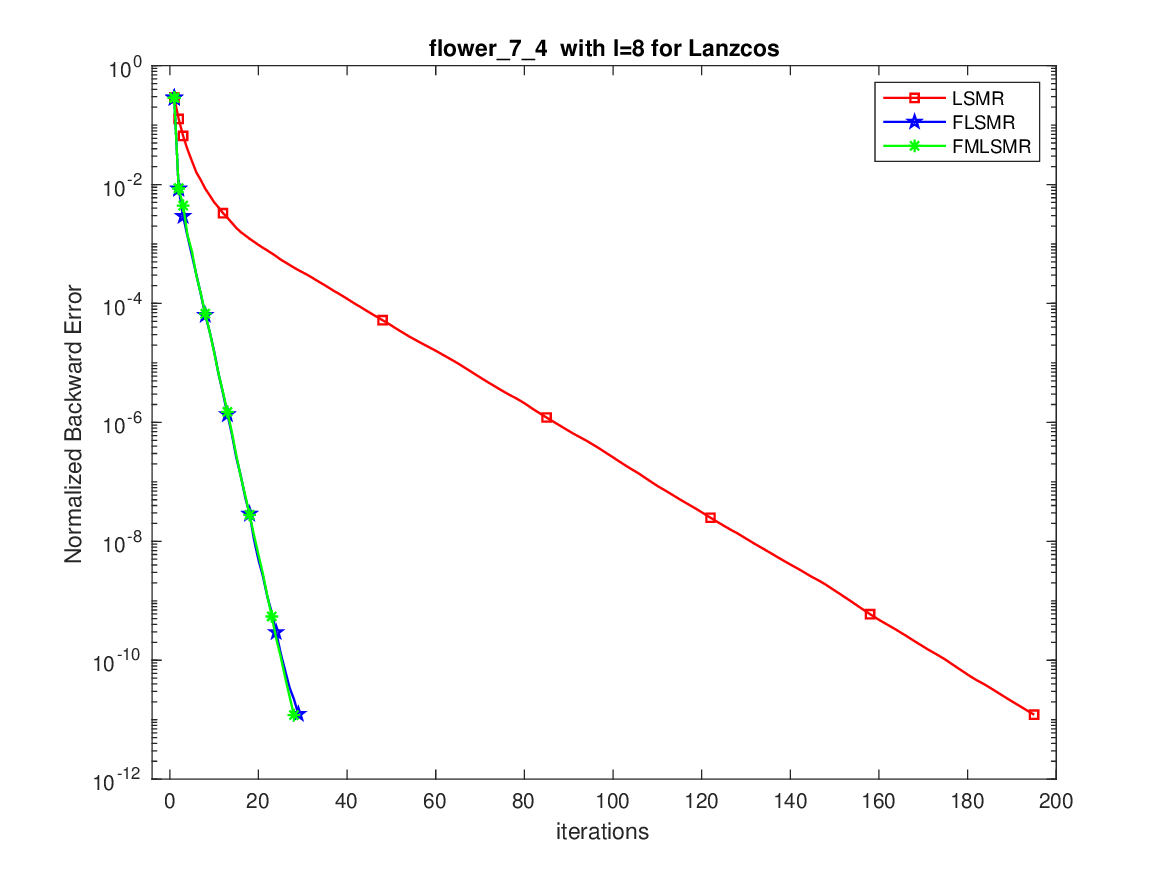}}\\
		\end{tabular}\par
	}\vspace{-0.15 cm}
	\caption{NRes ({\em left panel\/}) and normalized backward error $\|\hat{E}\|_2/\|A\|_1$ ({\em right panel\/})  for {\tt well1850}, {\tt cat\_ears\_3\_4},
		and {\tt flower\_7\_4}. }
	\label{fig:test41}
\end{figure}

\begin{figure}
	{\centering
		\begin{tabular}{cc}
			\hspace{-0.3 cm}
			\resizebox*{0.44\textwidth}{0.235\textheight}{\includegraphics{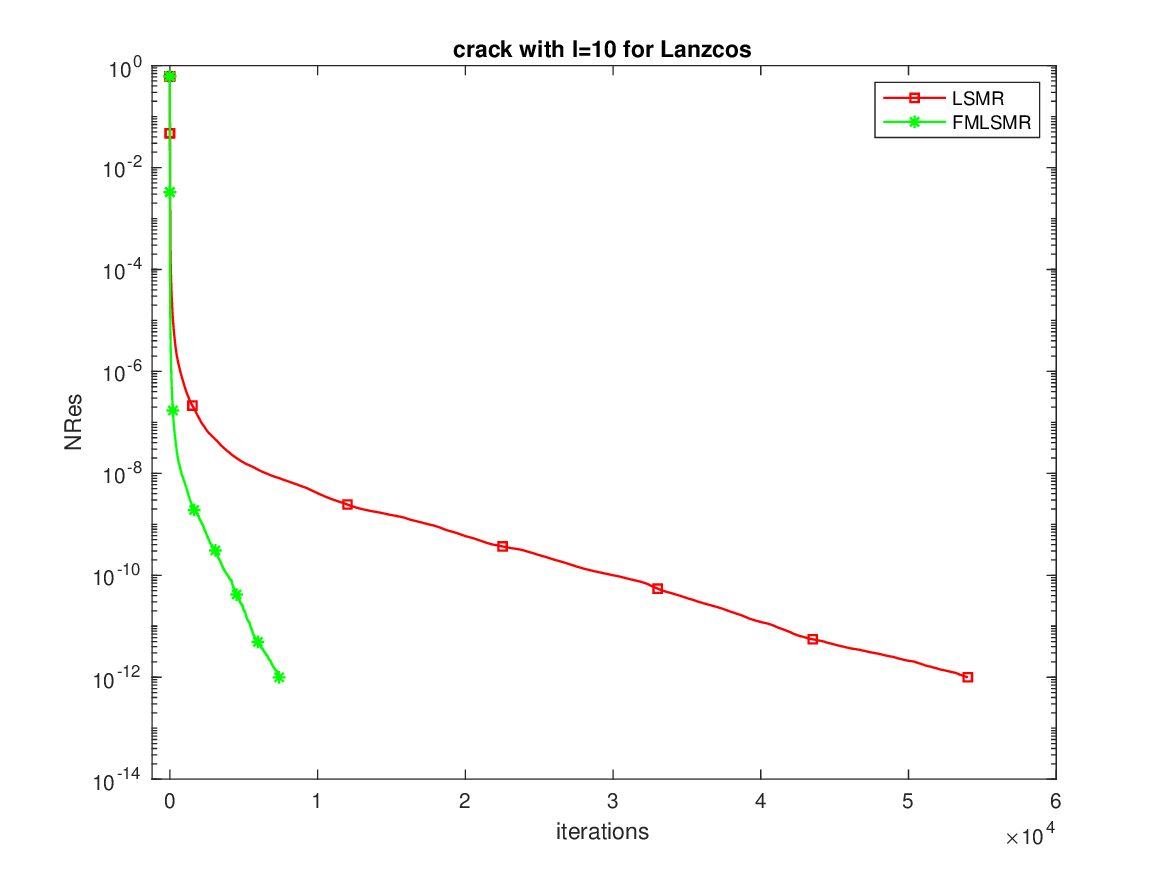}}
			&  \hspace{-0.5 cm}
			\resizebox*{0.44\textwidth}{0.235\textheight}{\includegraphics{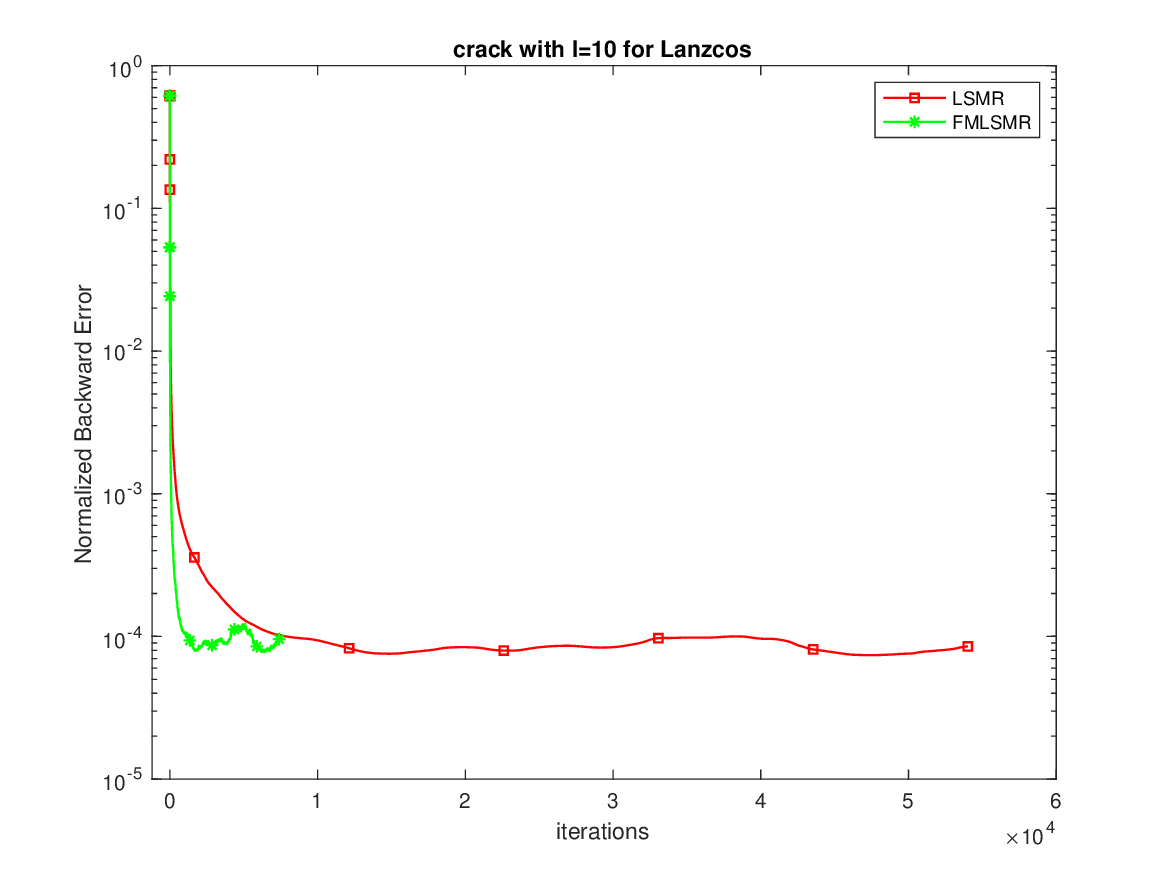}}\\
			\hspace{-0.3 cm}
			\resizebox*{0.44\textwidth}{0.235\textheight}{\includegraphics{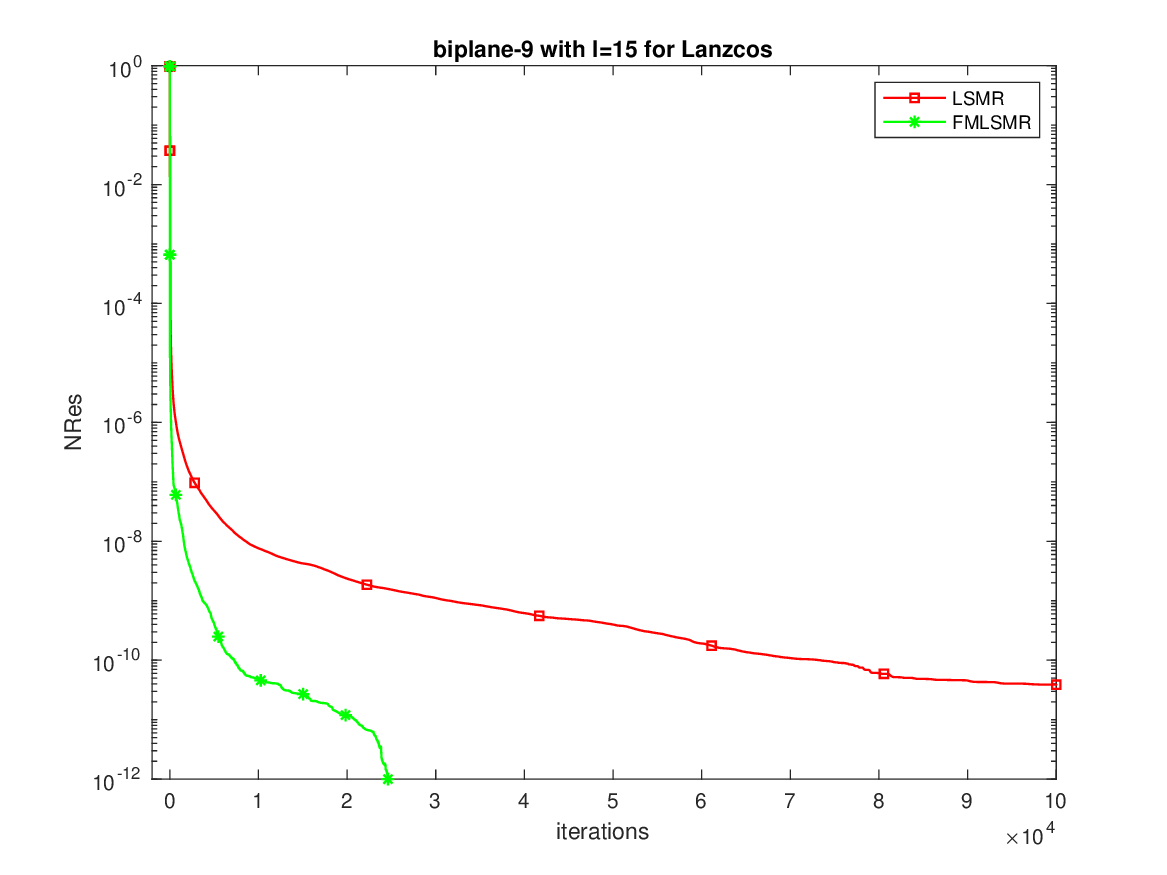}}
			&  \hspace{-0.5 cm}
			\resizebox*{0.44\textwidth}{0.235\textheight}{\includegraphics{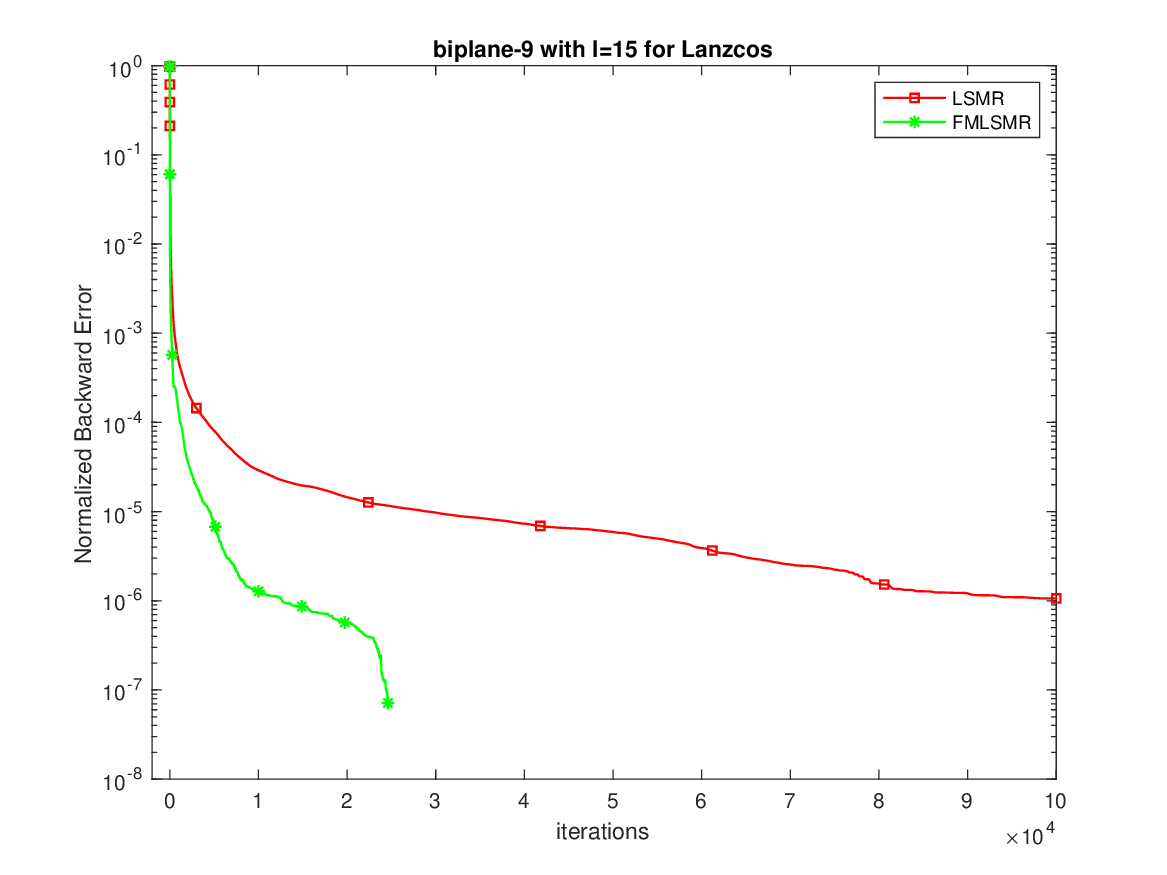}}\\
			\hspace{-0.3 cm}
			\resizebox*{0.44\textwidth}{0.235\textheight}{\includegraphics{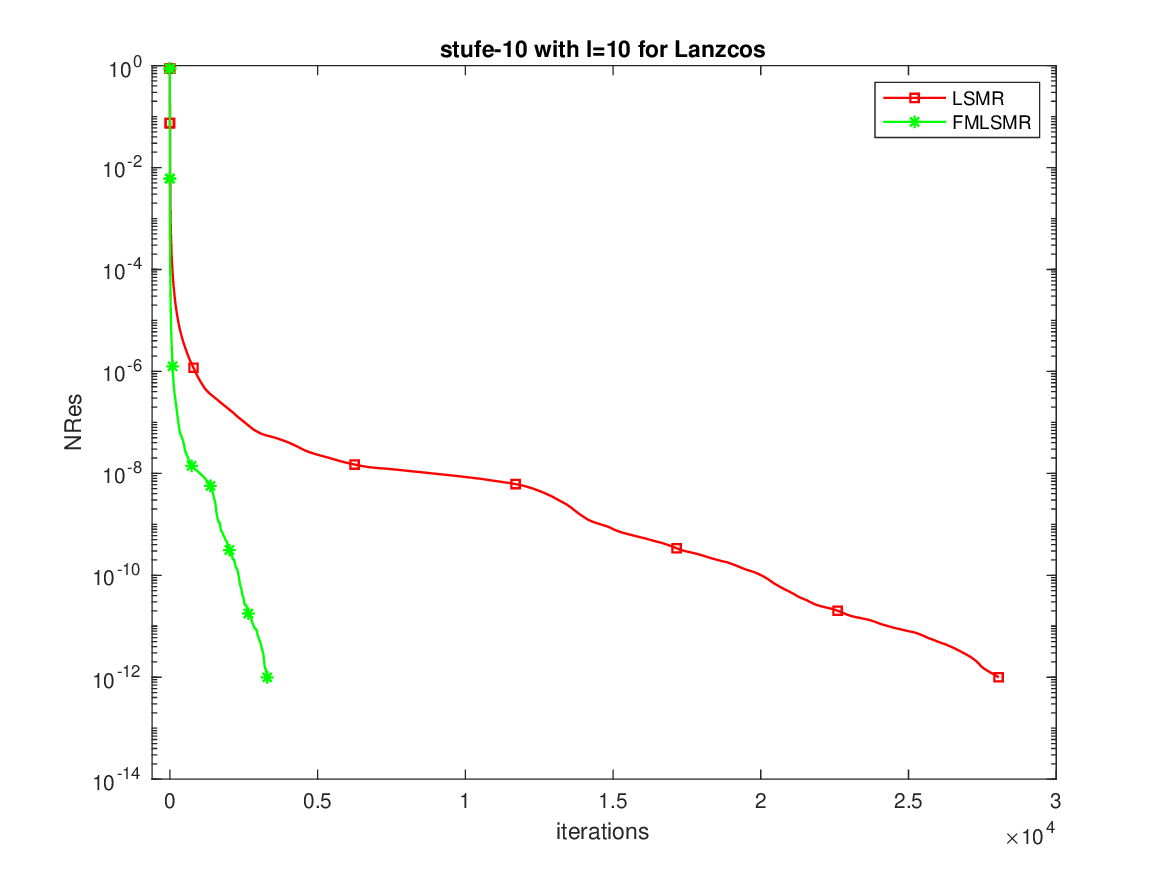}}
			&  \hspace{-0.5 cm}
			\resizebox*{0.44\textwidth}{0.235\textheight}{\includegraphics{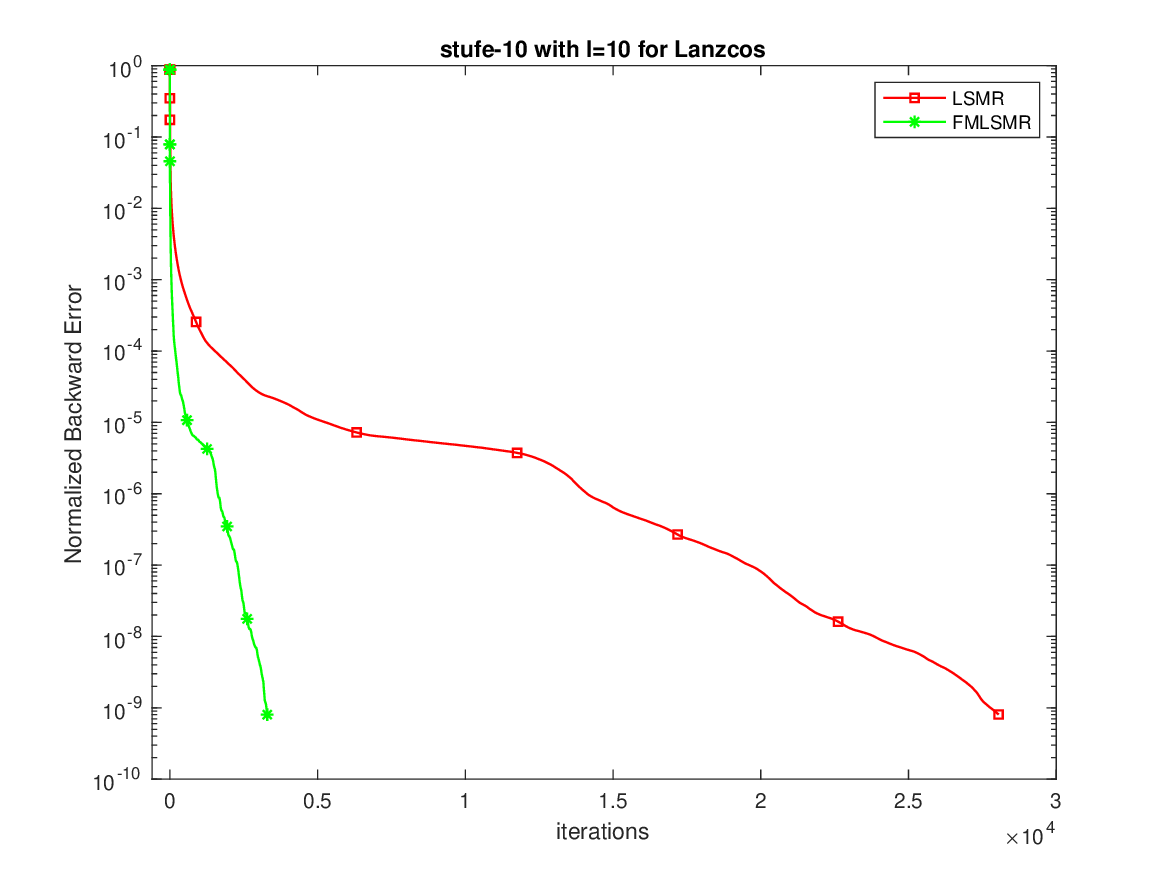}}\\
		\end{tabular}\par
	}\vspace{-0.15 cm}
	\caption{NRes ({\em left panel\/}) and normalized backward error  $\|\hat{E}\|_2/\|A\|_1$ ({\em right panel\/})  for {\tt crack}, {\tt biplane-9},
		and {\tt stufe-10}. }
	\label{fig:test42}
\end{figure}

Table~\ref{tbl:cycles4} collects the numbers of iterations by the methods on the eight problems, where the best results appear in boldface and for the places marked with
``--'' means that a method fails to solve a corresponding problem, i.e., satisfying \eqref{eq:NRes} within $10^5$
iterations.
Once again  the parameter $\ell$ in Table~\ref{tbl:cycles4} indicates that the $\ell$-step Lanczos method
is used to solve (\ref{eq:normal_type}) in FLSMR and FMLSMR.
%	\Red{
	%		\marginpar{\tiny ??}
	%		Since $\wtd {V}_k$ is not necessarily orthonormal in FMLSMR, we don't have to reorthogonalize $\wtd {V}_k$.
	%		For LSMR and FLSMR, the reorthogonalization is not applied, either.
	%	}
NRes and backward errors are displayed in Figures~\ref{fig:test41} and \ref{fig:test42} for six of the eight problems.
The total CPU times for each example are shown in Table~\ref{tbl:CPU}.
\begin{table}[h]
	\renewcommand{\arraystretch}{1.0}
	\caption{CPU Time (sec.)}\label{tbl:CPU}
	\centering
	\begin{tabular}{|c|c|c|c|c|c|}
		\hline
		ID &{matrix}                  &$\ell$ &LSMR      & FLSMR      &FMLSMR    \\ \hline\hline
		1  &{\tt well1850}         &8  & {\bf 0.0321}       & 0.4880    & 0.0518 \\
		2 &{\tt cat\_ears\_3\_4}  & 8 & 0.2999& 0.4761 & {\bf 0.1220}\\
		3  &\texttt{delaunay\_n16}   &30   & --   &  --  & {\bf1306.3191} \\
		4  &\texttt{biplane-9}   &15   & -- & --     & {\bf 452.9867}\\
		5  &\texttt{flower\_7\_4}  &8 & 0.4535   & 1.0227   &{\bf  0.4042}\\
		%				5 & \texttt{cs4}  &30  & 42528    & 1591   &2492\\
		%				6 & \texttt{cti} & 30 & 44621& 1665& 2542\\
		6 & \texttt{crack} &10 &  30.8903   &   --     & {\bf{28.7795}}    \\
		7 & \texttt{fe\_body}  &30 & -- & --   &    {\bf 517.2462}   \\
		8 & \texttt{stufe-10}  & 10 &  23.8165&  --& {\bf 23.0792}\\	
		
		\hline
	\end{tabular}
\end{table}

We have the following observations from Table~~\ref{tbl:cycles4} and \ref{tbl:CPU}, and Figures~\ref{fig:test41} and \ref{fig:test42}.	
\begin{itemize}
	\item The FMLSMR can solve all eight problems while LSMR succeeds on five of them and FLSMR on only three. Overall, FMLSMR has the best performance in terms of iteration numbers and CPU times, except for {\tt well11850} on which both 
	FLSMR and FMLSMR and fast while FLSMR holds an edge. 
	Specifically, for {\tt well11850} and {\tt stufe-10}, LSMR and FMLSMR has comparable performance in computational time.  For {\tt cat\_ears\_3\_4} and {\tt flower\_7\_4},
	FMLSMR is better than LSMR in terms of the number of iterations and CPU time.
	On {\tt delaunay\_n16}, {\tt biplane-9}, {\tt{fe\_body}}, {\tt{brack2}}, and {\tt stufe-10},
	FLSMR fails to satisfies the stopping criteria even for hours. According to Table~\ref{tbl:storage}, as the number of iterations increase, FLSMR uses much more storage and spends more on orthogonalization. The plots in the left column of Figures~\ref{fig:test41} and \ref{fig:test42} demonstrate a consistent decrease in relative residual across all methods. Notably, FMLSMR exhibits the fastest convergence among all. Therefore, considering the storage advantage of FMLSMR and its simple implementation, we can say that FMLSMR is a very good choice, especially for difficult problems, over FLSMR and LSMR.

	\item  The plots in the right column of Figures~\ref{fig:test41} and \ref{fig:test42} show
	backward errors \eqref{eq:BErr} for selective problems, and they display very similar patterns to that of NRes \eqref{eq:NRes}. Both $\|\hat{E}^{\rm FMLSMR}\|$ and $\|\hat{E}^{\rm FLSMR}\|$
	are less than $ \|\hat{E}^{\rm LSMR}\|$,  which indicates FMLSMR and FLSMR compute more accurate solutions than
	LSMR does for the same number of iterations.
	%		 Additionally, this result shows that FMLSMR is a stable numerical method.
	%  \item Figure \ref{fig:test43} shows the NRes $v..s$ flops for all examples. Together with \ref{tbl:CPU}, we find that FLSMR costs the most computational cost and the CPU time for all examples. The total cost and CPU time of FMLSMR are less than that of FLSMR. Considering the storage advantage of FMLSMR and its simple implementation, we see that FMLSMR is a good choice and it is more like a balance between FLSMR and LSMR.
\end{itemize}

\section{Conclusion}\label{sec:5}
In this paper, we present a new method, the Flexible Modified LSMR (FMLSMR), which integrates the key ideas from the Modified LSMR and Flexible GMRES algorithms. We conduct a theoretical analysis of the Modified LSMR and compare it with the Flexible LSMR (FLSMR) when using a given fixed preconditioner. Through numerical experiments, we illustrate the efficiency of FMLSMR from various angles. The advantages of our method in terms of storage and computational cost position it as a promising numerical method for tackling challenging problems in practical applications.

\bibliographystyle{plain}
\bibliography{references}

\end{document}